\documentstyle[twoside, 12pt]{article}
\renewcommand{\baselinestretch}{1.3}
\newtheorem{remark}{Remark}

\newtheorem{assumption}{Assumption}
\def\qed{ \ \vrule width.2cm height.2cm depth0cm\smallskip}
\def \ind{1\!\!1}

\newcommand{\brm}{\begin{rem}}
\newcommand{\ermq}{\end{rem}}

\newcommand{\ba}{\begin{array}}
\newcommand{\ea}{\end{array}}
\newcommand{\be}{\begin{equation}}
\newcommand{\ee}{\end{equation}}
\newcommand{\bea}{\begin{eqnarray}}
\newcommand{\eea}{\end{eqnarray}}
\newcommand{\beaa}{\begin{eqnarray*}}
\newcommand{\eeaa}{\end{eqnarray*}}

\def \R{I\!\!R}

\def\g{\gamma}
\def\d{\delta}
\def\e{\varepsilon}

\def\l{\lambda}

\def\t{\tau}
\def\f{\varphi}

\def\o{\omega}

\def\D{\Delta}

\def\cD{{\cal D}}

\def\cF{{\cal F}}

\def\cH{{\cal H}}

\def\cP{{\cal P}}

\def\cU{{\cal U}}

\def\no{\noindent}

\def\ms{\medskip}

\def\q{\quad}
\def\qq{\qquad}

\def\cds{\cdots}

\def\bF{{\bf F}}

\def\qed{ \hfill \vrule width.25cm height.25cm depth0cm\smallskip}
\newcommand{\dfnn}{\stackrel{\triangle}{=}}
\newcommand{\basa}{\begin{assumption}}
\newcommand{\easa}{\end{assumption}}

\newcommand{\bas}{\begin{assum}}
\newcommand{\eas}{\end{assum}}

\def\esssup{\mathop{\rm esssup}}

\def\cds{\cdots}

\def\dis{\displaystyle}

\def\bF{{\bf F}}

\newtheorem{thm}{Theorem}[section]
\newtheorem{lem}[thm]{Lemma}
\newtheorem{cor}[thm]{Corollary}
\newtheorem{prop}[thm]{Proposition}
\newtheorem{rem}[thm]{Remark}

\newtheorem{defn}[thm]{Definition}
\newtheorem{assum}[thm]{Assumption}

\title{The Starting and Stopping Problem under Knightian Uncertainty and Related Systems of Reflected BSDEs}
\author{Said
Hamad\`ene\thanks{Universit\'e du Maine, D\'epartement de
Math\'ematiques, Equipe Statistique et Processus, Avenue Olivier
Messiaen, 72085 Le Mans, Cedex 9, France. e-mail:
hamadene@univ-lemans.fr}\,\,\ and \,\,Jianfeng Zhang\thanks{USC
Department of Mathematics, S. Vermont Ave, KAP 108, Los Angeles, CA
90089, USA. e-mail:jianfenz@usc.edu. Research supported in part by NSF grants DMS
04-03575 and DMS 06-31366. Part of the work was done while this author was visiting Universit\'e du Maine, whose hospitality is greatly appreciated.}}
\begin{document}
\date{\today}
\maketitle

\begin{abstract}{\it
This article deals with the starting and stopping problem under
Knightian uncertainty, i.e., roughly speaking, when the probability
under which the future evolves is not exactly known. We show that
the lower price of a plant submitted to the decisions of starting and stopping
is given by a solution of a system of two reflected
backward stochastic differential equations (BSDEs for short). We
solve this latter system and we give the expression of the optimal
strategy. Further we consider a more general system of $m$ ($m\geq
2$) reflected BSDEs with interconnected obstacles. Once more we show
existence and uniqueness of the solution of that system.} $\qed$
\end{abstract}
{\bf AMS Classification subjects}: 60G40 ; 93E20 ; 62P20 ; 91B99.
\medskip

\no {$\bf Keywords$}: Real options ; Backward SDEs ; Reflected
BSDEs; Snell envelope; Stopping time ; Starting and stopping;
Switching. \ms

\no {\bf 0. Introduction:} \setcounter{equation}{0} We first
introduce through an example the standard starting and stopping (or
$switching$)  problem which has attracted a lot of interests during
the last decades (see the long list of bibliography and the
references therein).

Assume that a power plant produces electricity whose selling price,
as we know, fluctuates and depends on many factors such as consumer
demand, oil prices, weather and so on. It is also well known that
electricity cannot be stored and when produced it should be almost
immediately consumed. Therefore for obvious economic reasons,
electricity is produced only when there is enough profitability in
the market. Otherwise the power station is closed up to time when
the profitability is coming back, $i.e.$, till the time when the
market selling price of electricity reaches a level which makes the
production profitable again. Then for this power station there are
two modes, operating and closed. Accordingly, a management strategy
of the station is an increasing sequence of stopping times
$\delta=(\tau_n)_{n\geq 0}$ ($\tau_0=0$ and for any $n\geq 0$,
$\tau_n\leq \tau_{n+1}$). At time $\tau_n$, the manager switches the
mode of the station from its current one to the other. However
making a change of mode is not free and generates expenditures.

Suppose now that we have an adapted stochastic process
$X=(X_t)_{t\leq T}$ which stands for either the market electricity
price or factors which determine the price. 
When the
power station is run under a strategy $\delta=(\tau_n)_{n\geq 0}$,
its yield is given by a quantity denoted $J(\delta)$ which depends
also on $X$ and many other parameters such as utility functions,
expenditures, ... . Therefore the main problem is to find a
management strategy $\delta^*=(\tau_n^*)_{n\geq 1}$ such that for
any $\delta$ we have $J(\delta^*)\geq J(\delta)$, $i.e.$
$J(\delta^*)=\sup_\delta J(\delta)$. Once determined, the strategy
$\delta^*$ gives the optimal way of running the power plant and, as
a by-product, the real constant $J(\delta^*)$ is nothing else but
the fair price of the power plant in the energy market.

The two-mode starting and stopping problems attracted a lot of
research activity (see $e.g.$ \cite{[BO],[BO1],[BS], [DX], [DH],
[D], [DP],[DZ],[DZ2], [GP], hj, hib, [KMZ], porchet, porchet2, tri,
tri1, dz}, ... and the references therein).

Recently, Hamad\`ene and Jeanblanc \cite{hj} consider a finite
horizon two-modes when the price processes are only adapted to the
filtration generated by a Brownian motion. Porchet et $al.$ in
\cite{porchet} have considered the same problem with exponential
utilities and allow for the manager the possibility to invest in a
financial market. Djehiche and Hamad\`ene \cite{[DH]} studied also
this problem but the model integrates the risk of default of the
economic unit. Let us also mention the work by Hamad\`ene and Hdhiri
\cite{hib} where the set up of those latter papers is extended to
the case where the price processes of the underlying commodities are
adapted to a filtration generated by a Brownian motion and an
independent Poisson process.

Finally note that this two-mode switching problem models also
industries, like copper or aluminium mines,..., where parts of the
production process are temporarily reduced or shut down when e.g.\@
fuel, electricity or coal prices are too high to be profitable to
run them. A further area of applications includes Tolling
 Agreements (see Carmona and Ludkovski \cite{[CL]} and Deng and
Xia \cite{[DX]} for more details).

The natural extension of the two mode starting and stoping problem,
is the case where there are more than two modes for the production.
This problem has been recently considered by several authors amongst
we can quote Carmona and Ludkovski \cite{[CL]}, Djehiche et al.
\cite{dhp} and Porchet et $al.$ \cite{porchet2}.
\medskip

The studies quoted above, however, assume that future uncertainty is
characterized by a certain probability measure $P$ over the states
of nature. This turn out to assume that the firm is in a way
$certain$ that future market conditions are governed by this
particular probability measure $P$. The notion of Knightian
uncertainty introduced by F.H. Knight \cite{knight} assumes that it
is not granted that future uncertainty is characterized by a single
probability measure $P$ but other probabilities $P^u$, $u\in {\cal
U}$, are also likely. Usually those probabilities $P^u$ are supposed
not far from $P$. This notion will be defined later. Therefore one
of the main issues is, $e.g.$, related to the fair price of the
power plant in the market. If this latter quantity does not exist
what could be the lower price of the plant in accordance with the
sur-replication concepts well-known in mathematical finance.

To make things more clear suppose that the process $X$ is the price
of electricity in the energy market and assume that its dynamics is
given by the following standard differential equation:
$$
dX_t=X_t(r_tdt+\sigma_tdB_t),t\leq T \mbox{ and }X_0=x>0$$ where
$(B_t)_{t\leq T}$ is a Brownian motion, $r\dfnn(r_t)_{t\leq T}$ is
the spot interest rate and finally $(\sigma_t)_{t\leq T}$ the
volatility of the electricity price. So if the parameters $r$ and
$\sigma$ are known then the price of the power plant is just given
by
$\sup_\d J(\d)$. However usually it happens that the process $r$ is
not precisely known. We just have on it some confidence $i.e.$ we
know that $P-a.s.$, for any $t\in [0,T], r_t\in [-\kappa, \kappa]$
where $\kappa$ is a positive real constant which describes the
degree of Knightian uncertainty ($\kappa$-ignorance in the
terminology of Chen-Epstein (see \cite{chen})). Therefore possible
dynamics of the electricity price are the following:
$$
dX_t=X_t(u_tdt+\sigma_tdB_t),t\leq T \mbox{ and }X_0=x>0$$ where $B$
is once more a Brownian motion and $u\dfnn(u_t)_{t\leq T}$ is an
adapted stochastic process which takes its values in the compact set
$[-\kappa,\kappa]$ . In this case, things go on like incompleteness
in financial markets, we are just able to speak about the lower
price of the power plant which is given by the quantity:
\be\label{infsup} J^*\dfnn \sup_{\delta}\inf_{u}J(\delta,u),\ee
where $J(\delta,u)$ is the yield of the power plant when run under
the strategy $\delta$ and the future evolves according to the
probability $P^u$ for which $B$ is a Brownian motion. Mainly in this
work we aim at evaluating the quantity $J^*$ and providing a pair
$(\d^*,u^*)$ such that $J^*= J(\d^*,u^*).$

So in order to tackle our problem, using systems of reflected BSDEs
with $oblique$ $reflection$, we first provide a verification theorem
which shapes the problem under consideration. We show that when the
solution of the system exists it provides an optimal strategy
$(\d^*,u^*)$ of the switching problem under Knightian uncertainty.
Then we deal with a general system of $m$ ($m\geq 2$) reflected
BSDEs with oblique reflection for which we provide a solution. As a
by-product, we obtain that the verification theorem is satisfied and
therefore the switching problem solved. Further we address the
difficult issue of uniqueness of the solution of the general system.
Basically it turns out that the solution of that system can be
characterized as an optimal value for an appropriate switching
problem. Henceforth it is unique.

The idea of using reflected BSDEs in starting and stopping problems
with two modes appeared already in a previous work by Hamad\`ene
$\&$ Jeanblanc \cite{hj}. Then there were several works on this
subject using the same tool (see $e.g.$ \cite{[CL], porchet2}). In
\cite{[CL]}, the authors consider the multi-mode starting and
stopping problem. However they left open the question of the
existence of the solution of the system of reflected BSDEs with
oblique reflection, associated with the multi-state switching
problem. This question of existence/uniqueness is solved by Djehiche
et $al.$ in \cite{dhp}. Independent of our work, very recently Hu
$\&$ Tang \cite{ht} considered a quite more general, $w.r.t.$ the
one introduced in \cite{[CL]}, multi-dimensional reflected BSDE with
oblique reflection. They show existence and uniqueness of the
solution. However their framework is still somehow narrow since, due
to their techniques based on the use of local times and Tanaka's
formula, the assumptions they put on the data are rather stringent.

In this paper, using the notions of Snell envelope of processes
\cite{Elka,[H]} and the notion of smallest $g$-supermartingales
introduced by Mingyu $\&$ Peng \cite{xupeng} we provide new results,
$w.r.t.$ the ones of \cite{ht}, on existence/uniqueness of the
solution for the system of reflected BSDEs with oblique reflection.

This paper is organized as follows. In Section 1, we introduce the
problem and give some properties of the model. The quantities
$J(\d,u)$ are expressed by means of solutions of standard BSDEs
whose coefficients are not square integrable. Then we provide a
verification theorem which shapes the problem via systems of
reflected BSDEs with interconnected obstacles. The solution of the
system provides the pair $(\d^*,u^*)$ which achieves the $sup$ $inf$
in (\ref{infsup}). In Section 2, we consider a more general system
of reflected BSDEs, and show the existence of its solution. Finally
in Section 3 we characterize the solution as the optimal reward over
some appropriate set of strategies. This implies uniqueness of the
solution of the system. $\qed$

\section{The starting and stopping problem}

\subsection{The model}
Throughout this paper $(\Omega, {\cal F}, P)$ will be a fixed
complete probability space on which is defined a standard
$d$-dimensional Brownian motion $B=(B_t)_{0\leq t\leq T}$ whose
natural filtration is $(\cF_t^0\dfnn\sigma \{B_s, s\leq t\})_{0\leq
t\leq T}$. Let $ \bF\dfnn(\cF_t)_{0\leq t\leq T}$ be the completed
filtration of $(\cF_t^0)_{0\leq t\leq T}$ with the $P$-null sets of
${\cal F}$, hence $(\cF_t)_{0\leq t\leq T}$ satisfies the usual
conditions, $i.e.$, it is right continuous and complete.
Furthermore, let:
\begin{itemize}
\item[-] ${\cal P}$ be the $\sigma$-algebra on $[0,T]\times \Omega$ of
$\bF$-progressively measurable sets ;
\item[-] ${\cal H}^{p,l}$
be the set of $\cal P$-measurable and $\R^l$-valued processes
$\eta=(\eta_t)_{t\leq T}$ such that $E[\int_0^T|\eta_s|^pds]<\infty$
($p\geq 1$) ;
\item[-] ${\cal S}^2$  be the set of $\cal P$-measurable, continuous,
$R$-valued processes ${\eta}=({\eta}_t)_{t\leq T}$ such that
$E[\sup_{t\leq T}|{\eta}_t|^2]<\infty$ ; we denote by ${\cal A}$ the
subset of ${\cal S}^2$ which contains non-decreasing processes
$(K_t)_{\leq T}$ such that $K_0=0$;
\item[-] for any stopping time $\tau \in [0,T]$, ${\cal T}_\tau$ denotes
the set of all stopping times $\theta$ such that $\tau \leq \theta
\leq T$, $P-a.s.$
\item[-] the class [D] be the set of ${\cal P}$-measurable $rcll$ (right continuous with left
limits)  processes $V=(V_t)_{t\leq T}$ such that the set of random
variables $\{V_\t, \t \in {\cal T}_0\}$ is uniformly integrable.
\item[-] for any stopping time $\l$, $E_\l$ is the conditional
expectation with respect to ${\cal F}_\l$, $i.e.$, $E_\l[.]\dfnn
E_\l[. |{\cal F}_\l]$.
\end{itemize}

Let us now fix the data of the problem.
\medskip

$(i)$ Let $X\dfnn (X_t)_{0\le t\le T}$ be an $\cP$-measurable
process with values in $\R^k$  such that each component belongs to
${\cal S}^2$ (then $X$ is continuous). It stands for factors which
determine the market electricity price.
\medskip

$(ii)$  For $i=1,2$, let $\psi_i: (t,x)\in [0,T]\times \R^k \mapsto
\psi_i (t,x)\in \R$, be Borelean functions for which there exists a
constant $C$ such that $|\psi_i(t,x)|\leq C(1+|x|)$, $i=1,2$.
$\psi_1$ (resp. $\psi_2$) represents the utility function for the
power plant when it is in its operating (resp. close) mode. Actually
in a small interval $dt$, when the power plant is in its operating
(resp. closed) mode it generates a profit equal to $\psi_1(t,X_t)dt$
(resp. $\psi_2(t,X_t)dt$).
\medskip

$(iii)$ The switching of the power plant from one mode to another is
not free. Actually if at a stopping time $\tau$, the plant is
switched from the operating (resp. closed) mode to the closed (resp.
operating) one, the sunk cost is equal to $\f_1(\tau,X_\tau)$ (resp.
$\f_2(\tau,X_\tau)$) where the non-negative functions $\f_1,\f_2:
(t,x)\in [0,T]\times \R^k \mapsto \f_1(t,x),\f_2(t,x) \in \R^+$ are
continuous and linearly growing, $i.e.$, there exists a constant $C$
such that $|\f_i(t,x)|\leq C(1+|x|)$, $i=1,2$. Additionally they
verify $\f_1(t,x)+\f_2(t,x)>0$ for any $(t,x)\in [0,T]\times R^k$.
This latter requirement means that it is not free to make two
instantaneous switching at any time $t\leq T$.
\medskip

$(iv)$ Let $\delta=(\tau_n)_{n\geq 0}$ be an admissible management
strategy of the plant, $i.e.$, the $\tau_n$'s are $\bF$-stopping
times such that $\tau_n\leq \tau_{n+1}$ ($\tau_0=0$) for any $n\geq
0$ and $\lim_{n\rightarrow \infty}\tau_n=T$, P-$a.s.$. The set of
all admissible strategies will be denoted by $\cal D$. We assume
that the power plant is in its operating mode at the initial time
$t=0$. Therefore $\tau_{2n+1}$ (resp. $\tau_{2n}$) are the times
where the plant is switched from the operating (resp. closed) mode
to the closed (resp. operating) one.
\medskip

In the conventional model, $i.e.$, if we know that the future will
be governed by the probability measure $P$ the mean yield of the
power plant when run under the strategy $\delta=(\tau_n)_{n\geq 0}$
is given by :
$$
J(\delta)\dfnn E^P\Big\{\int_0^T\psi^\d(t,X_t)dt- A^\d_T\Big\},
$$
where $E^P$ is the expectation under the probability measure $P$,
\be
\label{psiAdelta}
\left\{\ba{lll}
\dis\psi^\d(t,x)\dfnn \sum_{n\ge 0}\Big[\psi_1(t,x)\ind_{[\t_{2n},\t_{2n+1})}(t)+\psi_2(t,x)\ind_{[\t_{2n+1},\t_{2n+2})}(t)\Big];\\
\dis A^\d_t \dfnn \sum_{n\geq
0}\Big[\f_1(\t_{2n+1},X_{\t_{2n+1}})\ind_{\{\tau_{2n+1}<t\}}+\f_2(\t_{2n+2},X_{\t_{2n+2}})\ind_{\{\tau_{2n+2}<t\}}\Big].
\ea\right.
\ee



\noindent Therefore the price of the power plant in the energy
market is just $\sup_{\delta \in {\cal D}}J(\delta)$.
\medskip

Knightian uncertainty amounts to suppose that we are not sure that
the future will evolve under the probability $P$ but other
probabilities $P^u$, $u\in {\cal U}$ (which we will precise later)
are also likewise. However we will suppose that those possible
probabilities $P^u$ are not far from $P$ in the sense that $P$ and
$P^u$ are equivalent. Actually we will assume that:
$$\frac{dP^u}{dP}=L^u_T\dfnn\exp\Big(\int_0^Tb(s,X.,u_s)dB_s-\frac{1}{2}\int_0^T|b(s,X.,u_s)|^2ds\Big)$$
where:

$(i)$ $u\dfnn (u_t)_{t\leq T}$ is an $\cP$-measurable process with
values in some compact set $U$. Hereafter $u$ will be called an
admissible control and the set of those controls will be denote by
${\cal U}$.

$(ii)$ $b: (t,x,u)\in [0,T]\times {\cal C}([0,T], \R^k)\times U
\mapsto b(t,x,u)\in \R^d$ is a Borel measurable and bounded
function. Moreover we assume that for any $(t,x)$, the mapping $u
\in U\mapsto b(t,x,u)\in \R^k$ is continuous and for any $u\in {\cal
U}$ the process $(b(t,X.,u_t))_{t\leq T}$ is ${\cal P}$-measurable.

Note that since the function $b$ is bounded then the random variable
$L^u_T$ has moment of any order, $i.e.$, for any $p\geq 1$,
$E[(L^u_T)^p]<\infty$ and if we set, for $t\leq T$,
$L_t\dfnn E[L^u_T|\cF_t]$ then the process $(L^u_t)_{t\leq T}$ satisfies
the following standard stochastic differential equation:
$$dL^u_t=L^u_tb(t,X,u_t)dB_t,t\leq T \,;\q L^u_0=1.$$

As previously mentioned, if the future evolves according to the
probability law $P^u$, $u\in {\cal U}$, then the fair price of the
power station in the energy market is given by: $$ J(u)=\sup_{\delta
\in {\cal D}}J(\delta,u)$$where \be \label{Jdu} J(\delta,u)\dfnn
E^u\Big\{\int_0^T\psi^\d(t,X_t)dt-A^\d_T\Big\},
\ee and $E^u$ is the expectation under $P^u$ and $\psi^\d, A^\d_T$
are defined by (\ref{psiAdelta}). However all the probability
measures are likewise therefore the selling lower price of the power
plant in the energy market is given by: \be \label{J*}
J^*\dfnn\sup_{\delta \in {\cal D}}J(\d);\q J(\d)\dfnn \inf_{u \in
{\cal U}}J(\delta,u). \ee Actually the quantity $J^*$ stands for the
optimal yield of the power plant in the worst case of evolution of
the future.  Therefore the problem we are interested in is to asses
the value $J^*$ and to find a pair $(\delta^*,u^*)$ such that
$$J^*=J(\d^*)=J(\delta^*,u^*)=\inf_{u \in \cU}J(\d^*,u).$$We note that, for any $u$, $J(\d^*, u)\ge
J^*$. However, for an arbitrary $\d$, in general we do not have
$J(\d, u)\ge J(\d, u^*)$. $\Box$
\begin{remark}

In the particular case where the process $X$ is the solution of the
following standard functional stochastic differential equation: \be
\label{sde1}dX_t=a(t,X.)dt+\sigma (t,X.)dB_t,\,\,t\leq T \mbox{ and
}X_0=x\ee with appropriate assumptions on the functions $a$ and
$\sigma$ in order to guarantee existence and uniqueness of the
solution of (\ref{sde1}), then thanks to Girsanov's Theorem we have:
$$dX_t=(a(t,X.)+\sigma (t,X.)b(t,X.))dt+\sigma (t,X.)dB^u_t,\,\,t\leq T \mbox{ and
}X_0=x$$where $B^u_t=B_t - \int_0^tb(s,X.,u_s)ds, t\leq T$, which is
well known that it is a Brownian motion under the probability
measure $P^u$. $\qed$
\end{remark}

\subsection{Properties of the model}

We are going to simplify the problem and to show that we can focus
only on a restricted set of strategies which satisfy appropriate
integrability conditions. So for any admissible strategy
$\delta=(\t_n)_{n\geq 0}$ $\in {\cal D}$ let us recall
(\ref{psiAdelta}), (\ref{Jdu}) and (\ref{J*}). Note that $\psi^\d,
A^\d$ do not depend on $u$ and $A^\d$ is $rcll$.

Now for $p\geq 1$ let us set
$$
{\cal D}_p\dfnn \{\delta \in {\cal
D},\mbox{ such that }\sup_{u\in {\cal U}}E^u[(A^\d_T)^p]<\infty\};\q
{\cal D}'\dfnn \cup_{p>1}{\cal D}_p.
$$
It follows that if $\delta \in {\cal D}-{\cal D}_1$, we have
$J(\delta)=-\infty$ since the process $(\psi^\d(t,X_t))_{t\leq T}$
belongs to $L^1(dt\times dP^u)$ due to the facts that
$(\psi_i(t,X_t))_{t\leq T}$,  $i=1,2$, belongs to ${\cal H}^{2,1}$
and that the random variable $L^u_T$ has moments of any order with
respect to the probability measure $P$. As a consequence, in our
objective to evaluate and characterize the quantity $J^*=\sup_{\d
\in {\cal D}}\inf_{u\in {\cal U}}J(\delta,u)$, we can discard the
admissible strategies $\delta$ which do not belong to ${\cal D}_1$.
\medskip

Next we introduce the Hamiltonian of the problem which is defined
by: for any $(t,x,u,z)\in [0,T]\times {\cal C}([0,T],\R^k)\times
U\times \R^d$,
$$
 H(t,x,u,z)\dfnn zb(t,x,u) \mbox{ and }H^*(t,x,z) \dfnn
\inf_{u\in U}H(t,x,u,z). $$ Since $b(t,x,u)$ is bounded then the
function $H$ and $H^*$ are uniformly Lipschitz w.r.t. $z$.
Additionally, thanks to Benes's selection Theorem, there exists a
measurable function $u^*:(t,x,z)\in [0,T]\times {\cal
C}([0,T],\R^k)\times \R^{d}\mapsto u^*(t,x,z) \in U$ such that:
$$
H^*(t,x,z)= [zb(t,x,u)]_{u=u^*(t,x,z)}.$$

We are now going to express the yields $J(\delta,u)$ by the means of
solutions of BSDEs whose coefficients are not square integrable.
Actually we have:
\begin{prop}\label{bsdeloc}
$(i)$ Let $\delta \in {\cal D}_1$ and $u\in {\cal U}$, then there
exists a unique pair of processes
$(Y^{\delta,u},Z^{\delta,u})$ such that the process
$(Y^{\delta,u}-A^\d)L^u$ is of class [D],
$\int_0^T|Z^{\d,u}_s|^2ds<\infty$  $a.s.$, and finally for any
$t\leq T$ we have:
$$
Y^{\delta,u}_t=\int_t^T(\psi^\d(s,X_s)+H(s,X_s,u_s,Z^{\delta,u}_s))ds-\int_t^TZ^{u,\delta}_sdB_s-(A^{\d}_T-A^\d_t)
$$
Moreover for any $t\leq T$ we have:
$$Y^{\delta,u}_t=E^u[\int_t^T\psi^\d(s,X_s)ds-(A^\d_T-A^\d_t)|\cF_t].$$

$(ii)$ For any $\delta \in {\cal D}'$, there exist $q>1$ and a
unique pair of processes $(Y^{\delta},Z^{\delta})$
such that: \be \label{eqq} \left\{
\begin{array}{l}
\dis E\Big\{\sup_{t\leq T}|Y^\d_t|^q +(\int_0^T|Z^\d_s|^2ds)^{q/2}\Big\}<\infty;\\
\dis Y^{\delta}_t=\int_t^T(\psi^\d(s,X_s)+H^*(s,X,Z^{\delta}_s))ds-\int_t^TZ^{\delta}_sdB_s-(A^{\d}_T-A^\d_t),t\leq
T.
\end{array}\right.
\ee Moreover for any $t\leq T$, $\dis Y^\d_t = {\rm essinf}_{u\in
{\cal U}} Y^{\d,u}_t$. In particular, $J(\d)=Y^\d_0$  and the optimal argument is
$(u^*(t,X,Z^\d_t))_{t\leq T}$.
\end{prop}
\noindent $Proof:$ $(i)$ Let $\d$ be a strategy which belongs to
${\cal D}_1$ and $u\in {\cal U}$.  Therefore we have
$E[L^u_TA^\d_T]=E^u[A^\d_T]<\infty$. Besides the process
$(L^u_t\psi^\d(t,X_t))_{t\leq T}$ belongs to $L^1(dt\otimes dP)$.
Henceforth thanks to the result by Briand et al. (\cite{briand},
Theorem 6.3, pp.18) related to solutions of BSDEs whose coefficients
belong only to $L^1$, there exists a unique pair of processes
$\tilde{Y}^{\delta,u}$ of class [D] and $\tilde{Z}^{\delta,u}$ such
that $E[(\int_0^T|\tilde{Z}^{\d,u}_s|^2ds)^\gamma]<\infty$, for any
$ \gamma \in ]0,1[$, which satisfy:
$$
\tilde{Y}^{\delta,u}_t=-L^u_TA^\d_T+\int_t^TL^u_s\psi^\d(s,X_s)ds
-\int_t^T\tilde{Z}^{u,\delta}_sdB_s, t\leq T.$$ Let us set now for
$t\leq T$,
$$Y^{\d,u}_t\dfnn \tilde{Y}^{\delta,u}_t(L^u_t)^{-1}+A^\d_t;\q
Z^{\d,u}_t\dfnn (L^u_t)^{-1}[ \tilde{Z}^{u,\d}_t-\tilde{Y}^{u,\d}_t
b(t,X,u_t)]. $$
 First note that $Y^{\d,u}$ is finite since
$A^\d_T<\infty$, P-$a.s.$ due to the equivalence of the probability
measures $P$ and $P^u$. Moreover $\int_0^T|Z^{\d,u}_s|^2ds <\infty$,
P-$a.s.$. Finally the process $(Y^{\delta,u}-A^\d)L^u$ is just
$\tilde{Y}^{\d,u}$ which belongs to class [D]. Using now It\^o's
formula for $Y^{\d,u}$ we get: $\forall t\leq T$,
$$Y^{\delta,u}_t=\int_t^T(\psi^\d(s,X_s)+H(s,X_s,u_s,Z^{\delta,u}_s))ds-\int_t^TZ^{u,\delta}_sdB_s-(A^{\d}_T-A^\d_t).$$
It remains to show that $Y^{\d,u}_t$ is just the conditional payoff
after $t$. Actually let $\l_n$ be the following stopping time:
$$\l_n\dfnn \inf\{t\geq 0, \int_0^t|Z^{\d,u}_s|^2ds \geq
n\}\wedge T.$$ Therefore
$$Y^{\d,u}_{t\wedge
{\l_n}}=E^u\Big\{Y^{\d,u}_{\l_n}-A^\d_{{\l_n}}+\int_{t\wedge
{\l_n}}^{\l_n}\psi^\d(s,X_s)ds +A^{\d}_{t\wedge {\l_n}}|\cF_{t\wedge
{\l_n}}\Big\}.$$ But the sequence of stopping times $(\l_n)_{n\geq
0}$ converges to $T$ and $L^u(Y^{\d,u}-A^\d)$ belongs to class [D],
therefore $Y^{\d,u}_{\l_n}-A^\d_{{\l_n}}\rightarrow -A^\d_T$ in
$L^1(dP^u)$. Besides the second term in the conditional expectation
converges also in $L^1(dP^u)$ to $\int_{t}^T \psi^\d(s,X_s)ds
+A^{\d}_{t}$. It follows that:
$$Y^{\delta,u}_t=E^u\Big\{\int_t^T\psi^\d(s,X_s)ds-(A^\d_T-A^\d_t)|\cF_t\Big\},\q \forall t\leq T,$$
which is the desired result.
\medskip

Let us now focus on $(ii)$. Let $\delta$ be a strategy of ${\cal
D}'$, therefore there exists $p>1$ such that $\sup_{u\in {\cal
U}}E^u[(A^\d_T)^p]<\infty$. As the moments of any order of
$(L^u_T)^{-1}$, $u\in {\cal U}$, exists then there exists $q>1$ such
that $E[(A^\d_T)^q]<\infty$. Now using once more the result by
Briand et al. (\cite{briand}, Theorem 4.2, pp.11) related to BSDEs
in $L^q$ $(q\in ]1,2[)$ there exists a pair of processes
$(\tilde{Y}^\d,{Z}^\d)$ such that:
$$\left\{
\begin{array}{l}
\dis E\Big\{\sup_{t\leq T}|\tilde Y^\d_t|^q +(\int_0^T|Z^\d_s|^2ds)^q\Big\}<\infty;\\
\dis \tilde{Y}^{\delta}_t=-A^{\d}_T+\int_t^T(\psi^\d(s,X_s)+
H^*(s,X_s,{Z}^{\delta}_s))ds-\int_t^T{Z}^{\delta}_sdB_s,\,t\leq T.
\end{array}\right.
$$
Now let us set $Y^\d=\tilde{Y}^\d+A^\d_t$, then the pair
$(Y^\d,Z^\d)$ is solution of the BSDE (\ref{eqq}).

Next for any $t\leq T$,
$H^*(t,X,Z^\d_t)=H(t,X,Z^\d_t,u^*(t,X,Z^\d_t))$ and since
$(Y^\d-A^\d)L^{u^*}$ belongs to class [D] (note that
$u^*=(u^*(t,X,Z^\d_t))_{t\leq T})$ then thanks to $(i)$ we have:
$$
Y^{\delta}_t=E^{u^*}\Big\{\int_t^T\psi^\d(s,X_s)ds-(A^\d_T-A^\d_t)|\cF_t\Big\}=Y^{\d,u^*}_t,\forall
t\leq T.$$ Next let $u\in {\cal U}$. Then for any $t\leq T$,
\beaa
Y^\d_t-Y^{\delta,u}_t&=& \int_t^T(H^*(s,X,Z^{\delta}_s)-H(s,X_s,u_s,Z^{\delta,u}_s))ds-\int_t^T(Z^{\delta}_s-Z^{\delta,u}_s)dB_s\\
&=&
\int_t^T(H^*(s,X,Z^{\delta}_s)-H(s,X_s,u_s,Z^{\delta}_s))ds-\q\int_t^T(Z^{\delta}_s-Z^{\delta,u}_s)dB^u_s.
\eeaa
 As $(Y^\d-Y^{\delta,u})L^u$ is of class [D] and since
$H^*(s,X,Z^{\delta}_s)-H(s,X_s,u_s,Z^{\delta}_s)\leq 0$ therefore,
arguing as previously by using appropriate stopping times, we obtain
$Y^\d_t-Y^{\delta,u}_t\leq 0$ for any $t\leq T$. Henceforth it holds
that:
$$\dis Y^\d_t = {\rm essinf}_{u\in {\cal U}} Y^{\d,u}_t, t\leq T,$$ and the
optimal argument is $u^*=(u^*(t,X,Z^\d_t))_{t\leq T}$. $\qed$
\medskip

We are now going to prove that the suprema of $J(\delta)$ over
${\cal D}_1$ and ${\cal D}'$ are the same. Actually we have:
\begin{prop}
\label{red}
$\dis\sup_{\delta \in {\cal D}_1}J(\delta)=\sup_{\delta \in {\cal
D}'}J(\delta).$
 \end{prop}$Proof$: For any $\d\in \cD_1$ and any
$n$, let $
\d^n \dfnn \{\t^n_i\}_{i\ge 0},
$ where
$$
\l_n \dfnn \inf\{t\geq 0: A^\d_t \ge n\}\wedge T;\q
\t^n_i \dfnn \left\{\ba{lll}
\t_i,\q {\rm if} ~~\t_i<\l_n;\\
T,\q {\rm if}~~\t_i\ge \l_n. \ea\right.
$$
It is obvious that the stopping times $\l_n \uparrow T$, and
$A^{\d^n}_T\le n$ and then $ \d^n \in \cD'.$

For any $u \in {\cal U}$,
$$
J(\d,u) =E^u\Big\{\int_0^T \psi^\d(t,X_t)dt - A^{\d}_T\Big\} \le
E^u\Big\{\int_0^T \psi^\d(t,X_t)dt - A^{\d^n}_T\Big\} \dfnn
J_n(\d,u).
$$
Note that
$$
\begin{array}{ll}
|J_n(\d,u)-J(\d^n,u)| &\le \dis E^u\Big\{\int_{\l^n}^T|\psi^\d(t,X_t)-
\psi^{\d^n}(t,X_t)|dt\Big\}\\
{}&\leq\dis
2\Big\{E^u[\int_0^T\max_{i=1,2}|\psi_i(t,X_t)|^pdt]\Big\}^{1/p}\Big\{E^u[(T-\l_n)]\Big\}^{1/q}\end{array}
$$  where $p\in ]1,2[$ and $q$ is its conjugate. But the right-hand side converges uniformly in $u\in \cU$
to $0$ as $n\rightarrow \infty$ since the processes
$(\psi_i(t,X_t))_{t\leq T}$ belong to $\cH^{2,1}$, $L^u_T$ have
moments of any order and $(b(t,X,u_t))_{t\leq T}$ is a uniformly
bounded process. Therefore we have:
$$
\lim_{n\to\infty}\sup_{u\in {\cal U}} |J_n(\d,u)-J(\d^n,u)| = 0.
$$
It follows that:
$$
J(\d,u) \le J_n(\d,u) = J_n(\d,u)-J(\d^n,u) + J(\d^n,u) \le
\sup_{u\in {\cal U}}|J_n(\d,u)-J(\d^n,u)| + J(\d^n,u).
$$
Minimizing now both hand-sides over $u\in \cU$, we get:
$$
J(\d) \le \sup_{u\in {\cal U}}|J_n(\d,u)-J(\d^n,u)| + J(\d^n) \le
\sup_{u\in {\cal U}}|J_n(\d,u)-J(\d^n,u)| + \sup_{\d\in\cD'}J(\d).
$$
Finally taking the limit as $n\to\infty$ to obtain the desired
result. $\qed$
\subsection{A verification theorem. Connection with reflected BSDEs}
 In order to tackle the problem which is
described in the previous part we are going to use the notion of
systems of backward stochastic differential equations with
reflecting barriers which we introduce now. \ms

Let us consider the following two dimensional reflected BSDEs:
 \be
 \label{RBSDE}
 \left\{\ba{lll}
 \dis Y^1,Y^2 \in {\cal S}^2, \,\,Z^1,Z^2 \in {\cal H}^{2,d} \mbox{ and }K^1,K^2 \in {\cal A},\\
 \dis Y^1_t = \int_t^T \Big[\psi_1(s,X_s)+ H^*(s,X_s,Z^1_s)\Big]ds - \int_t^T Z^{1}_s dB_s + K^1_T - K^1_t;\\
 \dis Y^2_t = \int_t^T \Big[\psi_2(s,X_s)+ H^*(s,X_s,Z^2_s)\Big]ds - \int_t^T Z^{2}_s dB_s + K^2_T - K^2_t;\\
 \dis Y^1_t \ge Y^2_t - \f_1(t,X_t);\q [Y^1_t - Y^2_t + \f_1(t,X_t)]dK^1_t = 0;\\
 \dis Y^2_t \ge Y^1_t - \f_2(t,X_t);\q [Y^2_t - Y^1_t + \f_2(t,X_t)]dK^2_t = 0.
 \ea\right.
 \ee
 For the moment we suppose that the processes $Y^i,Z^i,K^i$, $i=1,2$ exist. We leave the well-posedness and computation of (\ref{RBSDE}) to next
 section. Our main result of this section is the following theorem.

\begin{thm}
 \label{dopt}
 Assume $\f_1(t,x)+\f_2(t,x)>0$. Then $Y^1_0 = \sup_{\d \in \cD_1} \inf_{u\in \cU}J(\delta,u)$. Moreover, the optimal strategy
 $\d^*$ which belongs to  ${\cD}_1$ is given by
 $\t^*_0\dfnn 0$ and, for $n=0,\cds$,
   \beaa
   &&\t^*_{2n+1}\dfnn \inf\{t\ge \t^*_{2n}: Y^1_t = Y^2_t - \f_1(t,X_t)\}\wedge T;\\
   &&\t^*_{2n+2}\dfnn \inf\{t\ge \t^*_{2n+1}: Y^2_t = Y^1_t - \f_2(t,X_t)\}\wedge T.
   \eeaa
   \end{thm}
 {\it Proof.} First let us point out that thanks to Proposition \ref{red}, it is enough to show that
$Y^1_0 = \sup_{\d \in \cD'} \inf_{u\in \cU}J(\delta,u)$. So let
$\d=(\tau_n)_{n\geq 0} \in \cD'$ ($\tau_0=0$) and let us show that
we have $Y^1_0\ge Y^\d_0$. To this end, we
 define for $t\leq T$:
 \beaa
 \bar Y^\d_t &\dfnn&  \sum_{n=0}^\infty \Big[Y^1_t\ind_{[\t_{2n},\t_{2n+1})}(t) + Y^2_t\ind_{[\t_{2n+1},\t_{2n+2})}(t)\Big];\\
 \bar Z^\d_t &\dfnn&  \sum_{n=0}^\infty \Big[Z^1_t\ind_{[\t_{2n},\t_{2n+1})}(t) + Z^2_t\ind_{[\t_{2n+1},\t_{2n+2})}(t)\Big].
 \eeaa
Note that there is no problem of definition of the processes $\bar
Y^\d$ and $\bar Z^\d$ since the series are convergent (at least
pointwise). Besides $\bar Y^\d$ is $rcll$ and uniformly square
integrable and $\bar Z^\d$ belongs to $\cH^{2,d}$ for any admissible
strategy $\delta$. Moreover we have:
 \bea
 \dis \bar Y^\d_0 &=& Y^1_0 = Y^1_{\t_1} + \int_0^{\t_1} \Big[\psi_1(s,X_s)+ H^*(s,X_s,Z^1_s)\Big]ds - \int_0^{\t_1} Z^{1}_s dB_s + K^1_{\t_1}\nonumber\\
 \label{ineq1}
 \dis&\ge& Y^2_{\t_1} - \f_1({\t_1},X_{\t_1})\ind_{\{\t_1<T\}}\\
  \dis&&+ \int_0^{\t_1} \Big[\psi^\d(s,X_s)+H^*(s,X_s,\bar Z^\d_s)\Big]ds + \int_0^{\t_1} \bar Z^\d_s dB_s\nonumber\\
 \dis &=& Y^2_{\t_2} + \int_{\t_1}^{\t_2} \Big[\psi_2(s,X_s)+ H^*(s,X_s,Z^2_s)\Big]ds - \int_{\t_1}^{\t_2} Z^{2}_s dB_s + K^2_{\t_2}-K^2_{\t_1}\nonumber\\
 \dis&&  - \f_1({\t_1},X_{\t_1})\ind_{\{\t_1<T\}} + \int_0^{\t_1} \Big[\psi^\d(s,X_s)+ H^*(s,X_s,\bar Z^\d_s)\Big]ds - \int_0^{\t_1} \bar Z^\d_s dB_s\nonumber\\
 \label{ineq2}
 \dis&\ge& Y^1_{\t_2} - \f_2({\t_2},X_{\t_2})\ind_{\{\t_2<T\}}- \f_1({\t_1},X_{\t_1})\ind_{\{\t_1<T\}} \\
 \dis&&  + \int_0^{\t_2} \Big[\psi^\d(s,X_s)+ H^*(s,X_s,\bar Z^\d_s)\Big]ds - \int_0^{\t_2} \bar Z^\d_s dB_s.\nonumber
 \eea
 Repeat the  procedure as many times as necessary we get: for any
 $n\geq 0$,
\beaa
\bar Y^\d_0&\geq &Y^1_{\t_{2n+2}}-\sum_{k=0}^n\Big[
\f_1(\tau_{2k+1},X_{\tau_{2k+1}})\ind_{[\tau_{2k+1}<T]}+
 \f_2(\tau_{2k+2},X_{\tau_{2k+2}})\ind_{[\tau_{2k+2}<T]}\Big]\\
 &&+
 \int_0^{\t_{2n+2}} \Big[\psi^\d(s,X_s)+ H^*(s,X,\bar Z^\d_s)\Big]ds - \int_0^{\t_{2n+2}} \bar
 Z^\d_sdB_s
\eeaa
Taking now the limit as $n\rightarrow \infty$ and noting that $\t_n\uparrow T$, we obtain:
 $$
 \bar Y^\d_0 \ge \int_0^T \Big[\psi^\d(s,X_s)+ H^*(s,X,\bar Z^\d_s)\Big]ds - \int_0^T \bar Z^\d_s dB_s - A^\d_T.
 $$
 Following the same arguments we get, for any $t\leq T$,
 \be
 \label{barYd}
 \begin{array}{ll}
 \bar Y^\d_t& \ge \int_t^T \Big[\psi^\d(s,X_s)+ H^*(s,X,\bar Z^\d_s)\Big]ds - \int_t^T \bar Z^\d_s dB_s - [A^\d_T-A^\d_t]\\
 \end{array}
 \ee

 Here let us emphasize that up to now we did not use the fact
 that the strategy $\d$ belongs to $\cD'$ but only the fact that
 $\d$ is admissible. This remark will be useful later.

At this level we need $\d$ to be an element of $\cD'$. Actually let
us consider  the process $Y^\d$ defined in (\ref{eqq}). Then for any
$t\leq T$ we have, \be
 \label{barYd1}
 \begin{array}{ll}
 \bar Y^\d_t-Y^\d_t& \ge \dis \int_t^T \Big[H^*(s,X,\bar Z^\d_s)-H^*(s,X,Z^\d_s)\Big]ds - \int_t^T (\bar Z^\d_s-Z^\d_s) dB_s.\\
 {}&\geq \dis \int_t^T (\bar Z^\d_s-Z^\d_s) d\tilde B_s
 \end{array}
 \ee
where $\tilde B_.\dfnn B_.-\int_0^.\g_sds$ with
$$\g_s\dfnn\frac{H^*(s,X,\bar Z^\d_s)-H^*(s,X,Z^\d_s)}{\bar
Z^\d_s-Z^\d_s}\ind_{[\bar Z^\d_s-Z^\d_s\neq 0]}$$ which is a bounded
$\cP$-measurable process since the mapping $z\mapsto H^*(t,X,z)$ is
uniformly Lipschitz. Therefore, thanks to Girsanov's Theorem,
$\tilde B$ is a new Brownian motion under a new probability measure
$\tilde P$ equivalent to $P$ whose density w.r.t. $P$ is given by
$\tilde L$ which satisfies:$$ d\tilde L_t=\tilde L_t\g_tdB_t, \q
\tilde L_0=1.$$ Note that since the process $\g$ is bounded then the
random variable $\tilde L_T$ has moment of any order w.r.t. $P$.
Next we know that there exists a real constant $q>1$ such that
\\$E[(\int_0^T\{|\bar Z^\d_s|^2+|Z^\d_s|^2)ds)^{q/2}]<\infty$, then there exists another real constant $q'>1$ such that
$\tilde E[(\int_0^T\{|\bar Z^\d_s|^2+|Z^\d_s|^2)ds)^{q'/2}]<\infty$.
Therefore the stochastic integral $\int_0^.(\bar Z^\d_s-Z^\d_s)
d\tilde B_s$ is a actually a martingale. Going back now to
(\ref{barYd1}), taking expectation w.r.t. $\tilde P$ we obtain that
$\bar Y^\d_t-Y^\d_t\geq 0$ $\tilde P$-a.s. and then also $P$-a.s.
since the probabilities are equivalent. As this inequality is valid
for any $t\leq T$ and the processes $\bar Y^\d$ and $Y^\d$ are
$rcll$ then P-$a.s.$, for any $t\leq T$, $\bar Y^\d_t\geq
Y^\d_t={\rm essinf}_{u\in {\cU}}Y^{\d,u}_t$. \ms

It remains to prove $\d^*=(\t^*_n)_{n\ge 0}$ is optimal. First let us show that $\d^*$ is admissible,
$i.e.$, P-$a.s.$ $lim_{n\rightarrow \infty}\t^*_n=T$. Actually let
$\o$ be such that $lim_{n\rightarrow \infty}\t^*_n(\o)=\t^*(\o)<T$.
As the processes $Y^1,Y^2$, $(\varphi_1(t,X_t))_{t\le T}$ and
$(\varphi_2(t,X_t))_{t\le T}$ are continuous then for any $n\ge 0$
we have:
$$\begin{array}{l}Y^1_{\t^*_{2n+1}}(\o) = Y^2_{\t^*_{2n+1}}(\o) - \f_1(\t^*_{2n+1}(\o),X_{\t^*_{2n+1}}(\o))\mbox{ and
}\\ Y^2_{\t^*_{2n+2}}(\o) = Y^1_{\t^*_{2n+2}}(\o) -
\f_2(\t^*_{2n+2}(\o),X_{\t^*_{2n+1}}(\o)).\end{array}$$ We now let
$n$ tends to $+\infty$ and we obtain
$$\begin{array}{l}Y^1_{\t^*}(\o) = Y^2_{\t^*}(\o) - \f_1(\t^*(\o),X_{\t^*}(\o))\mbox{ and
} Y^2_{\t^*}(\o) = Y^1_{\t^*}(\o) -
\f_2(\t^*(\o),X_{\t^*}(\o))\end{array}$$ which obviously implies
that $\f_1(\t^*(\o),X_{\t^*}(\o))+\f_2(\t^*(\o),X_{\t^*}(\o))=0$
which is impossible. Therefore $P[\o: \lim_{n\rightarrow
\infty}\t^*_n(\o)<T]=0$ and the strategy $\d^*$ is admissible.

 On the other hand, 
 note that by definition $(Y^1,Y^2)$ are continuous processes, then
 $$
 Y^1_{\t^*_{1}} = Y^2_{\t^*_{1}} - \f_1(\t^*_{1},X_{\t^*_{1}})\ind_{\{\t^*_{1}<T\}};\q Y^2_{\t^*_{2}} = Y^1_{\t^*_{2}} - \f_1(\t^*_{2},X_{\t^*_{2}})\ind_{\{\t^*_{2}<T\}}.
 $$
 Moreover,
 $$
 K^1_{\t^*_{1}}=0;\q K^2_{\t^*_{2}}=K^2_{\t^*_{1}}.
 $$
 Therefore the inequalities (\ref{ineq1}) and  (\ref{ineq2}) become equalities. Following similar arguments and since $\d^*$ is admissible
 we have: for any $t\leq T$,
 $$
 \bar Y^{\d^*}_t = \int_t^T \Big[\psi^{\d^*}(s,X_s)- H(s,X_s,\bar Z^{\d^*}_s)\Big]ds - \int_t^T \bar Z^{\d^*}_s dB_s - [A^{\d^*}_T-A^{\d^*}_t].
 $$
Writing the equation for $t=0$ we deduce that
$E[(A^{\d^*}_T)^2]<\infty$ since $\bar Y^{\d^*}$ is uniformly square
integrable and $\bar Z^{\d^*}$ belongs to $\cH^{2,d}$. It follows
that there exists a constant $p\in ]1,2[$ such that $\sup_{u\in
\cU}E^u[(A^{\d^*}_T)^p]<\infty$ and then $\d^*$ belongs to ${\cal
D}'$. By the well-posedness of (\ref{eqq}) for elements of $\cD'$,
we get $\bar Y^{\d^*}_t = Y^{\d^*}_t$ since $\bar Y^{\d^*}$ and
$\bar Z^{\d^*}$ are adapted processes. In particular, $Y^1_0=\bar
Y^{\d^*}_0=Y^{\d^*}_0=\sup_{\d \in \cD'}Y^\d_0=\sup_{\d \in
\cD'}\inf_{u\in \cU}J(\d,u)=\sup_{\d \in \cD}\inf_{u\in
\cU}J(\d,u)=J^*$. Additionally $\d^*$ is optimal in $\cD_1$ since
$\cD' \subset \cD_1$. $\Box$
\begin{rem}: Thanks to Proposition \ref{bsdeloc}-(ii), the control $u^*=(u^*(t,X, \bar Z^\d))_{t\leq T}$ combined with the
strategy $\delta^*$ satisfy:
$$
Y^1_0=\bar Y^{\d^*}_0=Y^{\d^*}_0=J(\d^*,u^*)=\inf_{u\in
\cU}J(\d^*,u)=\sup_{\d \in \cD}\inf_{u\in \cU}J(\d,u).$$\qed
\end{rem}

\section{High Dimensional Reflected BSDEs: Existence}
\setcounter{equation}{0}
 As stated in Theorem \ref{dopt}, the solution of our original problem turns into solving the system
 of two reflected BSDEs (\ref{RBSDE}) whose obstacles are inter-connected and depend on the solution.
 Therefore in what follows we are going to deal with general systems of reflected BSDEs
 such that (\ref{RBSDE}) is just a particular case. Actually let us consider the following general system of
 RBSDEs: for $j=1,\cds,m$,
 \be
 \label{RBSDE0}
 \left\{\ba{lll}
 \dis Y^j\in {\cal S}^2, \,\,Z^j\in {\cal H}^{2,d} \mbox{ and }K^j\in {\cal A},\\
 \dis Y^j_t = \xi_j+\int_t^T f_j(s,Y^1_s,\cds,Y^m_s,Z^j_s)ds - \int_t^T Z^{j}_s dB_s + K^j_T - K^j_t;\\
 \dis Y^j_t \ge \max_{i\in A_j}h_{j,i}(t,Y^i_t);\q [Y^j_t - \max_{i\in A_j}h_{j,i}(t,Y^i_t)]dK^j_t =
 0;
 \ea\right.
 \ee
 where $A_j\subset \{1,\cds,m\}-\{j\}$, and the coefficients $f_j, h_{j,i}$ can depend upon $\o$. For simplicity we denote
 $\overrightarrow{Y_t}\dfnn (Y^1_t,\cds,Y^m_t)$, and similarly for other vectors. We emphasize that here $A_j$ can be
 empty and if so we take the convention that the maximum over the empty set, denoted as $\emptyset$, is $-\infty$.
 Then in this case $Y^j$ has no lower barrier and then we take $K^j=0$. Consequently, $Y^j$ satisfies the following BSDE
 without reflection:
 $$
 Y^j_t = \xi_j+\int_t^T f_j(s,\overrightarrow{Y}_s,Z^j_s)ds - \int_t^T Z^{j}_s dB_s,\,t\leq T.
 $$
 Also, for any $j$ we define
 \be
 \label{hjj}
 h_{j,j}(t,y) \dfnn y.
 \ee
 We note that the ${Y^j}$ of the solution of (\ref{RBSDE0}) satisfies
 \be
 \label{Yjj}
 Y^j_t \ge \max_{i\in A_j\cup\{j\}}h_{j,i}(t,Y^i_t).
 \ee

\begin{rem} The system we consider in (\ref{RBSDE0}) is
appropriate for multi-dimensional switching problems when from one
mode $j$ of the plant we are allowed to switch only to the modes
which belong to $A_j$. \qed
\end{rem}

Throughout this section we shall adopt the following assumption.
 \begin{assum}
 \label{assump} For any $j=1,\cds,m$, it holds that:

 (i) $\dis
 E\Big\{\int_0^T\sup_{\overrightarrow{y}: y_j=0}|f_j(t,\overrightarrow{y},0)|^2dt+|\xi_j|^2\Big\}<\infty.
 $

 (ii) $f_j(t,\overrightarrow{y},z)$ is uniformly Lipschitz continuous in $(y_j,z)$ and is continuous and increasing in $y_i$ for any $i\neq j$.

 (iii) For $i\in A_j$, $h_{j,i}(t,y)$ is continuous in $(t,y)$
 increasing in
 $y$, and $ h_{j,i}(t,y) \le y$. Moreover, if $j_2\in A_{j_1},\cds,j_k\in A_{j_{k-1}}, j_1\in A_{j_k}$,
 for any $y$, denote
 $$
 y_k\dfnn h_{j_k,j_1}(t,y),\q y_{k-1}\dfnn h_{j_{k-1},j_k}(t, y_{k}),\cds, y_1\dfnn h_{j_1,j_2}(t,y_2).
 $$
 Then we have
 \be
 \label{h}
 y_1 <y.
 \ee

 (iv) For any $j=1,...,m$, $\dis \xi_j\ge \max_{i\in A_j}h_{j,i}(T,\xi_i)$.\qed
 \end{assum}
 \begin{rem}The condition (\ref{h}) means that it is not free to make a circle of instantaneous
 switchings. It is satisfied if for example for any $i,j$, $h_{ij}(\omega,t,y)=y-c_{ij}(\o,t)$ with
 $c_{ij}(\o,t)>0$, $\forall t\leq T.$ $\Box$ \end{rem}
 Our main result is:
 \begin{thm}
 \label{wellposed0}
 Assume Assumption \ref{assump} holds true. Then RBSDE (\ref{RBSDE0})
 has at least one solution.
\end{thm}
$Proof:$ We use Picard iteration. First let us denote:
 \beaa
 \underbar f_j(t,y,z) \dfnn \inf_{\overrightarrow{y}:y_j=y}
 f_j(t,\overrightarrow{y},z) ~~\mbox{ and }~~
 \bar f_j(t,y,z) \dfnn \sup_{\overrightarrow{y}:y_j=y} f_j(t,\overrightarrow{y},z).
 \eeaa
 By Assumption \ref{assump} (i) and (ii), $ \underbar f_j, \bar f_j$ are uniformly Lipschitz continuous in $(y,z)$ and
 $$
 E\Big\{\int_0^T [|\underbar f_j(t,0,0)|^2 + |\bar
 f_j(t,0,0)|^2]dt\Big\}<\infty.
 $$
 Next, 
 let $(Y^{j,0},Z^{j,0})$ be the solution to the following BSDE without reflection:
  \be
  \label{Yj0}
 Y^{j,0}_t = \xi_j+\int_t^T \underbar f_j(s,Y^{j,0}_s,Z^{j,0}_s)ds - \int_t^T Z^{j,0}_s dB_s,\q j=1,\cds,m.
  \ee
  For $j=1,\dots, m$ and $n=1,2,\cds$, recursively define $Y^{j,n}$ via the following
  RBSDEs whose solution exits thanks to the result by El-Karoui et al. \cite{[EKal]}:
   \be
  \label{Yjn}
  \left\{\ba{lll}
 \dis Y^{j,n}_t = \xi_j- \int_t^T Z^{j,n}_s dB_s + K^{j,n}_T - K^{j,n}_t\\
  \dis\q+\int_t^T f_j(s,Y^{1,n-1}_s,\cds,Y^{j-1,n-1}_s,Y^{j,n}_s,Y^{j+1,n-1}_s\cds,Y^{m,n-1}_s,Z^{j,n}_s)ds;\\
 \dis Y^{j,n}_t \ge \max_{i\in A_j}h_{j,i}(t,Y^{i,n-1}_t);\q [Y^{j,n}_t - \max_{i\in A_j}h_{j,i}(t,Y^{i,n-1}_t)]dK^{j,n}_t = 0.
 \ea\right.
  \ee
  Note that, given $Y^{i,n-1}, i=1,\cds,m$, for each $j$ (\ref{Yjn}) is a one dimensional BSDE or reflected BSDE.
  Under Assumption \ref{assump}, (\ref{Yjn}) has a unique solution. Moreover, by comparison
  theorem (see e.g. \cite{[EKal]}, Theorem 4.1.) it is obvious that $Y^{j,1}\ge Y^{j,0}$. Then by induction one can easily show that $Y^{j,n}$ is increasing
  as $n$ increases.

 In order to obtain uniform estimates of $Y^{j,n}$, denote:
 $$
 \breve\xi\dfnn \sum_{j=1}^m |\xi_j|\mbox{ and } \breve f(t,y,z)\dfnn \sum_{j=1}^m |\bar f_j(t,y,z)|.
 $$
 Let $(\breve Y,\breve Z)$ be the solution to the following BSDE:
 $$
 \breve Y_t = \breve\xi + \int_t^T \breve f(s,\breve Y_s,\breve Z_s)ds -\int_t^T \breve Z_s dB_s.
 $$
 Denote, for $j=1,\cds,m$,
 $$
 \bar Y^j_t \dfnn \breve Y_t,\q \bar Z^j_t \dfnn \breve Z_t,\q \bar K^j_t\dfnn 0.
 $$
 Obviously $Y^{j,0}_t\le \bar Y^j_t$. Note that $(\bar Y^j,\bar Z^j,\bar K^j)$ satisfies
 $$
 \left\{\ba{lll}
 \dis \bar Y^{j}_t = \breve \xi + \int_t^T \breve f(s,\bar Y^j_s,\bar Z^j_s)- \int_t^T \bar Z^{j}_s dB_s + \bar K^{j}_T - \bar K^{j}_t;\\
 \dis \bar Y^{j}_t \ge \max_{i\in A_j}h_{j,i}(t,\bar Y^{i}_t);\q [\bar Y^{j}_t - \max_{i\in A_j}h_{j,i}(t,\bar Y^{i}_t)]d\bar K^{j}_t = 0.
 \ea\right.
 $$
 Once more apply the comparison theorem repeatedly, we get
 $$
 Y^{j,n}_t \le \breve Y_t,\q\forall n.
 $$
 Recall that $Y^{j,n}_t\ge Y^{j,0}_t$.  Then
 \be
 \label{barYest}
 \sum_{j=1}^m E\Big\{\sup_{0\le t\le T} |Y^{j,n}_t|^2\Big\} \le C<\infty,\q \forall n.
 \ee
 Moreover,
  $$
  E\Big\{\sup_{0\le t\le T} |[\max_{i\in A_j}h_{j,i}(t,Y^{i,n-1}_t)]^+|^2\Big\} \le E\Big\{\sup_{0\le t\le T} |[\max_{i\in A_j}Y^{i,n-1}_t]^+|^2\Big\} \le C.
  $$
  This further implies that
  \be
  \label{ZKjnest}
  E\Big\{\int_0^T|Z^{j,n}_t|^2dt + |K^{j,n}_T|^2\Big\} \le C,\q\forall j,n.
  \ee

 Now let $Y^j$ denote the limit of $Y^{j,n}$. By Peng's monotonic limit theorem \cite{peng} or \cite{xupeng},
  we know $Y^j$ is an $rcll$ process, and following similar arguments there one can easily show that
   there exist $(Z^j,K^j)$ such that
 \be
 \label{Yj}
 \left\{\ba{lll}
 \dis Y^j_t = \xi_j+\int_t^T f_j(s,\overrightarrow{Y}_s,Z^j_s)ds - \int_t^T Z^{j}_s dB_s + K^j_T - K^j_t;\\
 \dis Y^j_t \ge \max_{i\in A_j}h_{j,i}(t,Y^i_t).
 \ea\right.
 \ee
 Consider now the following RBSDEs whose solution exits thanks to the result by Hamad\`ene \cite{[H]} or Mingyu $\&$ Peng \cite{xupeng}:
 \be
 \label{Yjt-}
 \left\{\ba{lll}
 \dis \tilde Y^j_t = \xi_j- \int_t^T \tilde Z^{j}_s dB_s + \tilde K^j_T - \tilde K^j_t\\
 \dis\qq\qq+\int_t^T f_j(s,Y^1_s,\cds,Y^{j-1}_s,\tilde Y^j_s,Y^{j+1}_s,\cds,Y^m_s,\tilde Z^j_s)ds;\\
 \dis \tilde Y^j_t \ge \max_{i\in A_j}h_{j,i}(t,Y^i_t);\q [\tilde Y^j_{t-} - \max_{i\in A_j}h_{j,i}(t,Y^i_{t-})]d\tilde K^j_t = 0.
 \ea\right.
 \ee
 We note that (\ref{Yj}) and (\ref{Yjt-}) have the same lower barrier. Since  $\tilde Y^j_t$ is the smallest $f_j$-supermartingale
 with lower barrier $\dis \max_{i\in A_j}h_{j,i}(t,Y^i_t)$, we have $\tilde Y^j_t\le Y^{j}_t$ (see \cite{xupeng}, Theorem 2.1).
 On the other hand, since $ Y^{i,n-1}_t\le Y^i_t$ for any $(i,n-1)$, by the monotonicity of $h_{j,i}$ we get
 $$
 \max_{i\in A_j}h_{j,i}(t, Y^{i,n-1}_t) \le \max_{i\in A_j}h_{j,i}(t, Y^{i}_t).
  $$
 Then once more by comparison theorem for RBSDEs we have $Y^{j,n}_t\le \tilde Y^j_t$, which implies that
 $Y^j_t\le \tilde Y^{j}_t$. Therefore, $\tilde Y^j_t= Y^{j}_t$. This further implies that $dt\otimes dP$-$\tilde Z^j_t= Z^{j}_t$ and
 P-a.s. for any $t\leq T$, $\tilde K^j_t= K^{j}_t$, and that
 \be
 \label{Yjt}
 \left\{\ba{lll}
 \dis Y^j_t = \xi_j+\int_t^T f_j(s,\overrightarrow{Y}_s,Z^j_s)ds - \int_t^T Z^{j}_s dB_s + K^j_T - K^j_t;\\
 \dis Y^j_t \ge \max_{i\in A_j}h_{j,i}(t,Y^i_t),\q [Y^j_{t-} - \max_{i\in A_j}h_{j,i}(t,Y^i_{t-})]d K^j_t = 0..
 \ea\right.
 \ee

 Finally we show that $Y^j$ is continuous. We first note that, by  (\ref{Yjt}), $\D Y^j_t = -\D K^j_t \le 0$,
 and if $\D K^j_t\neq 0$, then $\dis Y^j_{t-} = \max_{i\in A_j}h_{j,i}(t,Y^i_{t-})$. It is obvious that $Y^j$
 is continuous when $A_j=\emptyset$. We now assume $\D Y^{j_1}_t\neq 0$ for some $j_1$ and $t$.
 Then $A_{j_1}\neq\emptyset$
 and $\D Y^{j_1}_t<0$. Note that in this case $\D K^{j_1}_t>0$, which further implies that
 $$
 Y^{j_1}_{t-} = \max_{i\in A_{j_1}}h_{j_1,i}(t,Y^i_{t-}).
 $$
 Let $j_2\in A_{j_1}$ be the optimal index, then
 $$
 h_{j_1,j_2}(t,Y^{j_2}_{t-}) = Y^{j_1}_{t-} > Y^{j_1}_{t}\ge \max_{i\in A_{j_1}}h_{j_1,i}(t,Y^i_{t})
 \ge h_{j_1,j_2}(t,Y^{j_2}_{t}).
 $$
 Thus $\D Y^{j_2}_t<0$, and therefore $A_{j_2}\neq\emptyset$. Repeat the arguments we obtain
 $j_k\in A_{j_{k-1}}$ and $\D Y^{j_k}_t<0$ for any $k$. Since each $j_k$ can take only values $1,\cds,m$, we may assume, without loss of generality that $j_1=j_{k+1}$ for some $k\ge 2$ (note again that $j_1\notin A_{j_1}$ and thus $j_2\neq j_1$). Then we have
 $$
 Y^{j_1}_{t-} = h_{j_1,j_2}(t,Y^{j_2}_{t-}), \cds,Y^{j_{k-1}}_{t-} = h_{j_{k-1},j_k}(t,Y^{j_k}_{t-}),\q Y^{j_{k}}_{t-} = h_{j_{k},j_1}(t,Y^{j_1}_{t-}).
 $$
 This contradicts with (\ref{h}).
 Therefore, all processes $Y^j$ are continuous.
\qed

By applying comparison theorem repeatedly, the following two results are direct consequence of Theorem \ref{wellposed0}, and their proofs are omitted.

\begin{cor}
\label{minsol}
The solution $\overrightarrow{Y}$ constructed in Theorem \ref{wellposed0} is the minimum solution to (\ref{RBSDE0}). That is, if $\overrightarrow{\tilde Y}$ is another solution to (\ref{RBSDE0}), then $Y^j_t\le \tilde Y^j_t, j=1,\cds,m$.
\end{cor}

\begin{cor}
\label{comparison0}
Assume $(\tilde\xi_j, \tilde f_j)$ also satisfy Assumption \ref{assump}, and
$$
f_j \le \tilde f_j,\q \xi_j\le \tilde\xi_j.
$$
Let $\overrightarrow{Y}$ and $\overrightarrow{\tilde Y}$ denote the
solution to (\ref{RBSDE0}) constructed in Theorem \ref{wellposed0},
with coefficients $(\xi_j, f_j, h_{j,i})$ and $(\tilde\xi_j,\tilde
f_j, \tilde h_{j,i})$, respectively. Then $Y^j_t\le \tilde Y^j_t,
j=1,\cds,m$.
\end{cor}

We now turn to the system (\ref{RBSDE}) and we have:
\begin{thm}
The system of reflected BSDEs (\ref{RBSDE}) has a unique solution.
\end{thm}
$Proof$: Existence is an immediate consequence of Theorem
\ref{wellposed0} through the properties satisfied by $\psi_1$,
$\psi_2$, $\varphi_1$, $\varphi_2$ and finally $H^*$ which make
Assumptions \ref{assump} fulfilled, especially the fact that
$\f_1(t,x)+\f_2(t,x)>0$ for any $(t,x)$. Uniqueness of $Y^1_0$ comes
from Theorem \ref{dopt}. Similarly one can prove the uniqueness of
$(Y^1_t,Y^2_t)$.
Uniqueness of $Z^1,Z^2$ is
a consequence of Doob-Meyer Decomposition, therefore we have
thoroughly uniqueness of $K^1$ and $K^2$. \qed

Another by-product of Theorem \ref{wellposed0} is that it provides
also existence of a solution of the system (\ref{RBSDE0}) considered
between two stopping times. This result is in particular useful to
show uniqueness of (\ref{RBSDE0}).
\ms

Actually let $\l_1$ and $\l_2$ be two stopping times such that
P-$a.s.$, $0\le \l_1 \le \l_2\le T$ and let us consider  the
following RBSDE over $[\l_1, \l_2]$: for $j=1,\cds,m$, P-$a.s.$, \be
 \label{RBSDEs}
 \left\{\ba{lll}
 \dis (Y^j_t)_{t\in [\l_1,\l_2]}\mbox{ continuous, }(K^j_t)_{t\in
 [\l_1,\l_2]} \mbox{ continuous and nondecreasing,} \\
  \dis K^j_{\l_1}=0, \mbox{ and }E\Big\{\sup_{t\in
 [\l_1,\l_2]}|Y^j_t|^2 +\int_{\l_1}^{\l_2}|Z^j_s|^2ds+(K^j_{\l_2})^2\Big\}<\infty;\\
 \dis Y^j_t = \xi^j_{\l_2}+\int_t^{\l_2} f_j(s,\overrightarrow{Y}_s,Z^j_s)ds -
 \int_t^{\l_2} Z^{j}_s dB_s + K^j_{\l_2} - K^j_t,\,\,\forall t\in [\l_1,\l_2]\,;\\
 \dis Y^j_t \ge \max_{i\in A_j}h_{j,i}(t,Y^i_t)\mbox{ and } [Y^j_t - \max_{i\in A_j}h_{j,i}(t,Y^i_t)]dK^j_t
 = 0,\,\forall t\in [\l_1,\l_2].
 \ea\right.
 \ee
 Then we have:
  \begin{thm}:
 \label{wellposed}
  Assume Assumption \ref{assump} holds true and that for $j=1,...,m$, $\xi^j_{\l_2}\in\cF_{\l_2}$ and satisfies:
  \be
  \label{terminal}
  E\{|\xi^j_{\l_2}|^2\} <\infty \mbox{ and }\,\xi^j_{\l_2} \ge
  \max_{i\in A_j} h_{j,i}(\l_2, \xi^i_{\l_2}).
  \ee
  Then the RBSDE (\ref{RBSDEs}) has a solution.
   \qed
  \end{thm}
\section{Uniqueness}
\setcounter{equation}{0} We now focus on uniqueness of the solution
of RBSDE (\ref{RBSDEs}), hence that of RBSDE (\ref{RBSDE0}). To do that we need a stronger assumption.
 \begin{assum}
 \label{assump2}
 (i) $f_j$ is uniformly Lipschitz continuous in all $y_i$.

 (ii)  If $i\in A_j, k\in A_i$, then $k\in A_j\cup\{j\}$. Moreover,
 \be
 \label{h2}
 h_{j,i}(t, h_{i,k}(t,y)) < h_{j,k}(t,y).
 \ee

 (iii) For any $i\in A_j$,
 \be
 \label{hLipschitz}
 |h_{j,i}(t,y_1)-h_{j,i}(t,y_2)|\le |y_1-y_2|.
 \ee


 \end{assum}
Note that these assumptions are satisfied if $A_j=\{1,\dots,
m\}-\{j\}$ for any $j=1,...,m$ and $h_{ij}(\o,t,y)=y-c_{ij}(\o,t)$
with $c_{ij}(\o,t)>0$ for any $t\leq T$, P-$a.s.$
\begin{thm}{\bf (Uniqueness)}
\label {uniqueness}

  (i) Assume Assumptions \ref{assump} and \ref{assump2} are in force. Then the solution to (\ref{RBSDEs}) is unique.

  (ii) Moreover, assume for $j=1,\dots,m$, $\tilde f_j$ satisfies Assumptions \ref{assump} and \ref{assump2},
  and $\tilde \xi^j_{\l_2}$ satisfies (\ref{terminal}). Let $(\tilde Y^j,\tilde Z^j)$  be the solution to RBSDE
  (\ref{RBSDEs}) corresponding to $(\tilde f_j,\tilde \xi^j_{\l_2})$. For
  $j=1,...,m$, denote,
  \be
  \label{DYxi}
 \D Y^j_t \dfnn Y^j_t - \tilde Y^j_t, \q \D \xi^j_{\l_2} \dfnn \xi^j_{\l_2}-\tilde \xi^j_{\l_2},\q
  \|\D f_t\|\dfnn \sum_{j=1}^m\sup_{(\overrightarrow{y},
  z)}|[f_j-\tilde f_j](t,\overrightarrow{y},z)|.
  \ee
 Then there exists a constant $C$, which is independent of $\l_1, \l_2$, such
 that:
 \be
 \label{DY}
 \max_{1\le j  \le m} |\D Y^j_{\l_1}|^2 \le E_{\l_1}\Big\{e^{C(\l_2-\l_1)}\max_{1\le j\le m}|\D
 \xi^j_{\l_2}|^2 + C \int_{\l_1}^{\l_2}\|\D
 f_t\|^2dt\Big\}.
 \ee
\end{thm}
\ms

The proof will be obtained after intermediary results. However
basically it uses an induction argument and a characterization of
$Y^j$ as a supremum over strategies $\delta$ of some processes
$Y^{j,\d}$ which are uniquely defined.

So assume Assumptions \ref{assump} and \ref{assump2} hold.
 Let $\mu$ denote the number of nonempty sets $A_j$ in (\ref{RBSDEs}), that is, the number of reflections
 in (\ref{RBSDEs}). We proceed by induction on $\mu$. First, when $\mu=0$, (\ref{RBSDEs}) becomes an $m$-dimensional BSDE without reflection.
 By standard arguments one can easily show that Theorem \ref{uniqueness} holds true. Now assume it is true for $\mu =m_1-1$
 for some $1\le m_1\le m$. For $\mu=m_1$, let $(Y^j,Z^j,K^j)$ be an arbitrary solution to (\ref{RBSDEs}).

 \subsection{Admissible strategies}
 We want to extend the arguments in Theorem \ref{dopt} to this case. The idea is to express $Y^j_t$ as
 the supremum of $Y^{j,\d}_t$, where $\d$ is an {\it admissible strategy} which we are going to define soon,
 and $Y^{j,\d}$ is the solution to a system of RBSDEs with $m_1-1$ reflections. Thus by induction
 $Y^{j,\d}$ is unique for each $(j,\d)$ and therefore $Y^j$ is unique.

 To motivate the definition of admissible strategy, we heuristically discuss how to find the
 ``optimal strategy", an analogue of the $\t^*_n$ in Theorem \ref{dopt}. A rigorous and more
 detailed argument will be given in \S3.3.

 Let $\t^*_0\dfnn \l_1$, and without loss of generality assume $A_1\neq\emptyset$. Set
 $$
 \t^*_1 \dfnn \inf\{t\ge \t^*_0: Y^1_t = \max_{i\in A_1} h_{1,i}(t, Y^i_t)\}\wedge \l_2.
 $$
 When $\t^*_1<\l_2$, we have
 $$
 Y^1_{\t_1^*} = \max_{i\in A_1} h_{1,i}(\t_1^*, Y^i_{\t_1^*}).
 $$
 That is, there exists an index, denoted as $\eta_1\in A_1$, such that
 $$
 Y^1_{\t_1^*} = h_{1,\eta_1}(\t_1^*, Y^{\eta_1}_{\t_1^*}).
 $$
 So, besides the stopping time $\t_1^*$, we need to keep track of the ``optimal index" $\eta_1$.
 At this point, let us denote $\eta_0\dfnn 1$. Note that, over $[\t_0^*, \t_1^*]$, it holds
 that:
 $$
 \left\{\ba{lll}
 \dis Y^{j}_t = Y^{j}_{\t_1^*}+\int_t^{\t_1^*} f_j(s,\overrightarrow{Y}_s,Z^{j}_s)ds- \int_t^{\t_1^*} Z^{j}_s dB_s + K^{j}_{\t_1^*} - K^{j}_t,~ j\neq \eta_0;\\
 \dis Y^{j}_t \ge \max_{k\in A_j}h_{j,k}(t,Y^{k}_t);\q [Y^{j}_t - \max_{k\in A_j}h_{j,k}(t,Y^{k}_t)]dK^{k}_t = 0,~ j\neq \eta_0;\\
 \dis Y^{\eta_0}_t =  Y^{\eta_0}_{\t_1^*}+\int_t^{\t_1^*} f_{\eta_0}(s,\overrightarrow{Y}_s,Z^{\eta_0}_s)ds- \int_t^{\t_1^*} Z^{\eta_0}_s dB_s.
 \ea\right.
 $$
 This is a system with only $m_1-1$ reflections.

 Now for $(\t_1^*,\eta_1)$, we need to consider two different cases.

 \noindent {\it \underline{ Case 1}.} Assume $A_{\eta_1}\neq\emptyset$. Then by considering $Y^{\eta_1}$
 over $[\t^*_1,\l_2]$ instead of $Y^{\eta_0}$ over $[\t_0^*,\l_2]$, similarly one can define $\t_2^*$ and
  $\eta_2\in A_{\eta_1}$, and see that $\overrightarrow{Y}$ satisfies a system with $m_1-1$ reflections
   over $[\t_1^*,\t_2^*]$, where the $\eta_1$-th equation has no reflection.

  \noindent {\it \underline{ Case 2}.} Assume $A_{\eta_1}=\emptyset$. In this case, the $\eta_1$-th equation
   has no reflection. Note that $Y^{\eta_0}_{\t_1^*} = h_{\eta_0,\eta_1}(\t_1^*, Y^{\eta_1}_{\t_1^*})$.
   Choose $\t_2^*$ ``close" to $\t_1^*$, then for any $t\in [\t_1^*,\t_2^*]$, we have
   $Y^{\eta_0}_t \approx h_{\eta_0,\eta_1}(\t_1^*, Y^{\eta_1}_t)$. On the other hand, by (\ref{h2})
   one can see that $Y^j_{\t_1^*} > h_{j,\eta_0}(\t_1^*, Y^{\eta_0}_{\t_1^*})$ for any $j$ such that
   $\eta_0\in A_j$. Since $\t_2^*$ is close to $\t_1^*$, let us assume  $Y^j_{t} > h_{j,\eta_0}(\t_1^*, Y^{\eta_0}_t)$
   for $t\in [\t_1^*,\t_2^*]$. So approximately, over $[\t_1^*,\t_2^*]$, $Y^j, j\neq\eta_0$ satisfy
  $$
 \left\{\ba{lll}
 \dis Y^{j}_t \approx Y^{j}_{\t_2^*}+\int_t^{\t_2^*} f_j(s,h_{1,\eta_1}(\t_1^*,Y^{\eta_1}_s), Y^2_s,\cds,Y^m_s,Z^{j}_s)ds- \int_t^{\t_2^*} Z^{j}_s dB_s + K^{j}_{\t_2^*} - K^{j}_t;\\
 \dis Y^{j}_t \ge \max_{k\in A_j-\{\eta_0\}}h_{j,k}(t,Y^{k}_t);\q [Y^{j}_t - \max_{k\in A_j-\{\eta_0\}}h_{j,k}(t,Y^{k}_t)]dK^{k}_t = 0.
 \ea\right.
 $$
 This is a system of $m-1$ equations with $m_1-1$ reflections, where we remove the equation for $Y^{\eta_0}$ completely.

 In order to move forward, we need to define $\eta_2$ so that $A_{\eta_2}\neq\emptyset$.
 It turns out that the best way is to set $\eta_2\dfnn \eta_0$. Then we can continue the procedure.

 Based on the above argument, let us introduce the following:
 \begin{defn}
 \label{admissible}
  $\d=(\t_0,\cds,\t_n; \eta_0,\cds,\eta_n)$ is called an admissible strategy if

 (i)  $\l_1=\t_0\le\cds\le\t_n\le\l_2$ is a sequence of stopping times;

 (ii) $\eta_0,\cds,\eta_n$ are random index taking value in $\{1,\cds,m\}$ such that $\eta_i\in\cF_{\t_i}$;

  (iii) $A_{\eta_0}\neq\emptyset$;

 (iv) If $A_{\eta_i}\neq\emptyset$, then  $\eta_{i+1}\in A_{\eta_i}$;

   (v) $A_{\eta_i}=\emptyset$, then $\eta_{i+1}\dfnn \eta_{i-1}$.
  \end{defn}

  \begin{rem}
  \label{etai+-1}
  By Definition \ref{admissible} (iii), $A_{\eta_i}=\emptyset$ implies that $i\ge 1$.
  Then (v) makes sense. Moreover,  in this case
  $A_{\eta_{i+1}}=A_{\eta_{i-1}}\neq \emptyset$.
  \end{rem}
 \subsection{Construction of $Y^\d$}
  For an admissible strategy $\d$, we construct $(Y^{\d,j}, Z^{\d,j})$ as follows.
   First, for $t\in [\t_n,\l_2]$ and $j=1,\cds,m$, set
    \be
    \label{Ydn}
    Y^{\d,j}_t\dfnn Y^{0,j}_t,\q Z^{\d,j}_t\dfnn Z^{0,j}_t,
    \ee
     where $(Y^{0,j},Z^{0,j})$ is the solution to (\ref{RBSDEs}) constructed
    in \S2. Then in particular we have
    \be
    \label{Y>hn}
    Y^{\d,j}_{\t_n} \ge \max_{i\in A_j}h_{j,i}(\t_n,
    Y^{\d,i}_{\t_n}),\q j=1,\cds,m.
    \ee

  For $i=n-1,\cds,0$, assume we have constructed $Y^{\d,j}_{\t_{i+1}-}$ for $j=1,\cds,m$,
  which we will do later. Note that $Y^{\d,j}$ may be discontinuous at $\t_{i+1}$.
   We define $(Y^{\d,j},Z^{\d,j})$ over $[\t_i, \t_{i+1})$ in two cases.
   \ms

  \noindent {\it \underline{ Case 1}.} If $A_{\eta_i}\neq\emptyset$, assume,
  \be
 \label{terminal1}
 Y^{\d,j}_{\t_{i+1}-} \ge \max_{k\in A_j} h_{j, k}(\t_{i+1}, Y^{\d,k}_{\t_{i+1}-}),\q  j\neq\eta_i.
 \ee
  We consider the following RBSDE by removing the constraint of the $\eta_i$-th equation:
 \be
 \label{Ydeta1}
 \left\{\ba{lll}
 \dis Y^{\d,j}_t = Y^{\d,j}_{\t_{i+1}-}+\int_t^{\t_{i+1}} f_j(s,\overrightarrow{Y}^{\d}_s,Z^{\d,j}_s)ds- \int_t^{\t_{i+1}} Z^{\d,j}_s dB_s + K^{\d,j}_{\t_{i+1}} - K^{\d,j}_t, ~j\neq \eta_i;\\
 \dis Y^{\d,j}_t \ge \max_{k\in A_j}h_{j,k}(t,Y^{\d,k}_t);~ [Y^{\d,j}_t - \max_{k\in A_j}h_{j,k}(t,Y^{\d,k}_t)]dK^{\d,k}_t = 0,~j\neq \eta_i;\\
 \dis Y^{\d,\eta_i}_t = Y^{\d,\eta_i}_{\t_{i+1}-}+\int_t^{\t_{i+1}} f_{\eta_i}(s,\overrightarrow{Y}^{\d}_s,Z^{\d,\eta_i}_s)ds- \int_t^{\t_{i+1}} Z^{\d,\eta_i}_s dB_s.
 \ea\right.
 \ee
 It is obvious that the $f_j, h_{j,i}, A_j$ here satisfy Assumptions \ref{assump} and \ref{assump2}.
  Since (\ref{Ydeta1}) has only $m_1-1$ reflections, by induction (\ref{Ydeta1}) has a unique solution $(Y^{\d, j}, Z^{\d, j}), j=1,\cds,m$ over $[\t_i,\t_{i+1})$.
  $\Box$
\ms

 \noindent {\it \underline{ Case 2}.} If $A_{\eta_i}=\emptyset$,
   by Remark \ref{etai+-1} we have $i\ge 1$ and $A_{\eta_{i-1}}\neq\emptyset$. Assume
  \be
 \label{terminal2}
 Y^{\d,j}_{\t_{i+1}-} \ge \max_{k\in A_j-\{\eta_{i-1}\}} h_{j, k}(\t_{i+1}, Y^{\d,k}_{\t_{i+1}-}), \q j\neq\eta_{i-1}.
 \ee
 We now omit the $\eta_{i-1}$-th equation and consider the following $m-1$ dimensional RBSDE with at most
$m_1-1$ reflections: for $j\neq \eta_{i-1}$,
 \be
 \label{Ydeta2}
 \left\{\ba{lll}
 \dis Y^{\d,j}_t = Y^{\d,j}_{\t_{i+1}-}- \int_t^{\t_{i+1}} Z^{\d,j}_s dB_s + K^{\d,j}_{\t_{i+1}} - K^{\d,j}_t\\
 \dis \q +\int_t^{\t_{i+1}} \tilde f_j(s,Y^{\d,1}_s,\cds,Y^{\d,\eta_{i-1}-1}_s,Y^{\d,\eta_{i+1}-1}_s,\cds,Y^{\d,m}_s,Z^{\d,j}_s)ds;\\
 \dis Y^{\d,j}_t \ge \max_{k\in A_j-\{\eta_{i-1}\}}h_{j,k}(t,Y^{\d,k}_t),\q [Y^{\d,j}_t - \max_{k\in A_j-\{\eta_{i-1}\}}h_{j,k}(t,Y^{\d,k}_t)]dK^{\d,k}_t = 0.
 \ea\right.
 \ee
 Here:
 \bea
 \label{tildef}
 &&\tilde f_j(t,y_1,\cds,y_{\eta_{i-1}-1},y_{\eta_{i-1}+1},\cds, y_n,z)\\
 &&\dfnn f_j(t,y_1,\cds,y_{\eta_{i-1}-1},h_{\eta_{i-1},\eta_i}(\t_i,y_{\eta_i}), y_{\eta_{i-1}+1},\cds, y_n,z).\nonumber
 \eea
 One can easily check that $\tilde f_j, h_{j,i}, A_j-\{\eta_{i-1}\}$ here satisfy Assumptions \ref{assump} and \ref{assump2}.
      Since (\ref{Ydeta2}) has at most $m_1-1$ reflections, by induction (\ref{Ydeta2}) has a unique solution
      $(Y^{\d, j}, Z^{\d, j}), j\neq \eta_{i-1},$ over $[\t_i,\t_{i+1})$.$\Box$
 \ms

 It remains to construct $Y^{\d,j}_{\t_{i+1}-}$ satisfying (\ref{terminal1}) or (\ref{terminal2}).  First, if $i+1=n$,
 set $Y^{\d, j}_{\t_{i+1}-} \dfnn Y^{0,j}_{\t_n}$ ; and if $\t_{i+1}=\l_2$, set $Y^{\d,j}_{\t_{i+1}-}\dfnn \xi^j_{\l_2}$. By (\ref{Ydn}) and (\ref{terminal}) we know both (\ref{terminal1}) and (\ref{terminal2}) hold true.
 Now assume $i<n-1$ and  $\t_{i+1}<\l_2$. Assume we have solved either (\ref{Ydeta1}) or (\ref{Ydeta2}) over $[\t_{i+1},\t_{i+2})$.
 \ms

\noindent {\it \underline{ Case 2}.} Assume $A_{\eta_i}=\emptyset$.
By Remark \ref{etai+-1} we know $i\ge 1, \eta_{i+1}=\eta_{i-1}$, and
$A_{\eta_{i+1}}\neq \emptyset$. Then we obtain $Y^{\d,j}_{\t_{i+1}}$
from (\ref{Ydeta1})  over $[\t_{i+1},\t_{i+2})$ satisfying:
 \be
 \label{Y>h3}
 Y^{\d,j}_{\t_{i+1}}\ge \max_{k\in A_j}h_{j,k}(\t_{i+1},Y^{\d,k}_{\t_{i+1}}),\q j\neq \eta_{i+1}=\eta_{i-1}.
 \ee
 Define
 \be
 \label{Yi-3}
 Y^{\d,j}_{\t_{i+1}-} \dfnn Y^{\d, j}_{\t_{i+1}},~~ j\neq \eta_{i-1}.
 \ee
 Then (\ref{terminal2}) follows immediately from (\ref{Y>h3}).
 $\Box$
 \ms

\noindent {\it \underline{ Case 1}.} Assume $A_{\eta_i}\neq\emptyset$. We
further discuss two cases. \ms

\noindent {\it \underline{ Case 1.1}.}  Assume
$A_{\eta_{i+1}}=\emptyset$. Then we obtain $Y^{\d,j}_{\t_{i+1}},
j\neq\eta_i$ from (\ref{Ydeta2}) over $[\t_{i+1},\t_{i+2})$
satisfying
 \be
 \label{Y>h2}
 Y^{\d,j}_{\t_{i+1}} \ge \max_{k\in A_j-\{\eta_{i}\}}h_{j,k}(\t_{i+1},Y^{\d,k}_{\t_{i+1}}),\q j\neq\eta_i.
 \ee
 Define
 \be
 \label{Yi-2}
Y^{\d, j}_{\t_{i+1}-} \dfnn Y^{\d, j}_{\t_{i+1}},~~ j\neq \eta_i;\q
  Y^{\d,\eta_{i}}_{\t_{i+1}-}\dfnn h_{\eta_i,\eta_{i+1}}(\t_{i+1}, Y^{\d,\eta_{i+1}}_{\t_{i+1}}).
 \ee
 By (\ref{Y>h2}), to prove (\ref{terminal1}) it suffices to show that
 \be
 \label{Y>h4}
 Y^{\d, j}_{\t_{i+1}} \ge h_{j,\eta_i}(\t_{i+1},h_{\eta_i,\eta_{i+1}}(\t_{i+1}, Y^{\d,\eta_{i+1}}_{\t_{i+1}})),\q \mbox{ if} ~~ \eta_i\in A_j.
 \ee
 By (\ref{h2}), we have
 $$
 h_{j,\eta_i}(\t_{i+1},h_{\eta_i,\eta_{i+1}}(\t_{i+1}, Y^{\d,\eta_{i+1}}_{\t_{i+1}})) < h_{j,\eta_{i+1}}(\t_{i+1}, Y^{\d,\eta_{i+1}}_{\t_{i+1}}).
 $$
 When $\eta_i\in A_j$, by Assumption \ref{assump2} (ii), we have $\eta_{i+1}\in [A_j-\{\eta_i\}]\cup\{j\}$. If $\eta_{i+1}\in A_j-\{\eta_i\}$, then (\ref{Y>h4}) follows (\ref{Y>h2}). If $\eta_{i+1}=j$, then (\ref{Y>h4}) follows (\ref{hjj}). So in both cases (\ref{Y>h4}) holds true, then so does (\ref{terminal1}).
 $\Box$

\noindent {\it \underline{ Case 1.2}.} Assume
$A_{\eta_{i+1}}\neq\emptyset$. Then we obtain $Y^{\d,j}_{\t_{i+1}}$ from
(\ref{Ydeta1}) over $[\t_{i+1},\t_{i+2})$ satisfying:
 \be
 \label{Y>h1}
 Y^{\d,j}_{\t_{i+1}}\ge \max_{k\in A_j}h_{j,k}(\t_{i+1},Y^{\d,k}_{\t_{i+1}}),\q j\neq \eta_{i+1}.
 \ee
 Define
 \be
 \label{Yi-1}
 \left.\ba{lll}
\dis Y^{\d, j}_{\t_{i+1}-} \dfnn Y^{\d, j}_{\t_{i+1}},\q j\neq \eta_i, \eta_{i+1};\\
 \dis Y^{\d,\eta_{i+1}}_{\t_{i+1}-} \dfnn Y^{\d,\eta_{i+1}}_{\t_{i+1}}\vee \max_{k\in A_{\eta_{i+1}}-\{\eta_i\}} h_{\eta_{i+1},k}(\t_{i+1}, Y^{\d,k}_{\t_{i+1}});\\
  \dis Y^{\d,\eta_{i}}_{\t_{i+1}-}\dfnn h_{\eta_i,\eta_{i+1}}(\t_{i+1}, Y^{\d,\eta_{i+1}}_{\t_{i+1}-}).
 \ea\right.
 \ee
 We now check (\ref{terminal1}) for $j\neq \eta_i$. First, for $j=\eta_{i+1}$, by (\ref{Yi-1}),
 $$
 Y^{\d,\eta_{i+1}}_{\t_{i+1}-} \ge \max_{k\in A_{\eta_{i+1}}-\{\eta_i\}} h_{\eta_{i+1},k}(\t_{i+1}, Y^{\d,k}_{\t_{i+1}-}).
 $$
 Moreover, if $\eta_i\in A_{\eta_{i+1}}$, by (\ref{h2}) and (\ref{hjj}) we have
 $$
 h_{\eta_{i+1},\eta_i}(\t_{i+1}, Y^{\d,\eta_i}_{\t_{i+1}-})= h_{\eta_{i+1},\eta_i}(\t_{i+1},h_{\eta_i,\eta_{i+1}}(\t_{i+1}, Y^{\d,\eta_{i+1}}_{\t_{i+1}-}))<Y^{\d,\eta_{i+1}}_{\t_{i+1}-}.
 $$
 So (\ref{terminal1}) holds true for $j=\eta_{i+1}$.

 Next, assume $j\neq\eta_i, \eta_{i+1}$, by (\ref{Y>h1}) and the first line in (\ref{Yi-1}) we have
 \be
 \label{Y>h11}
 Y^{\d,j}_{\t_{i+1}-}\ge \max_{k\in A_j-\{\eta_i,\eta_{i+1}\}}h_{j,k}(\t_{i+1},Y^{\d,k}_{\t_{i+1}-}).
 \ee
 If $\eta_{i+1}\in A_j$, recall the definition of $Y^{\d,\eta_{i+1}}_{\t_{i+1}-}$ in (\ref{Yi-1}). First, by (\ref{Y>h1}) we have
 $$
 h_{j,\eta_{i+1}}(\t_{i+1},Y^{\d,\eta_{i+1}}_{\t_{i+1}}) \le Y^{\d,j}_{\t_{i+1}}=Y^{\d,j}_{\t_{i+1}-}.
 $$
 Second, for any $k\in A_{\eta_{i+1}}-\{\eta_i\}$, similar to (\ref{Y>h4}) one can easily prove
 $$
 h_{j,\eta_{i+1}}(\t_{i+1},h_{\eta_{i+1},k}(\t_{i+1},Y^{\d,k}_{\t_{i+1}}))\le Y^{\d,j}_{\t_{i+1}}=Y^{\d,j}_{\t_{i+1}-}.
 $$
 Thus
\be
 \label{Y>h12}
 h_{j,\eta_{i+1}}(\t_{i+1},Y^{\d,\eta_{i+1}}_{\t_{i+1}-})\le Y^{\d,j}_{\t_{i+1}-}.
 \ee
 Finally, if $\eta_i\in A_j$, since $\eta_{i+1}\in A_{\eta_i}$, by Assumption \ref{assump2} (ii) we have $\eta_{i+1}\in A_j\bigcup\{j\}$.
 Then by (\ref{Y>h12}) and (\ref{h}) we have
  $$
  h_{j,\eta_{i}}(\t_{i+1},Y^{\d,\eta_{i}}_{\t_{i+1}-})
  =h_{j,\eta_{i}}(\t_{i+1},h_{\eta_i,\eta_{i+1}}(\t_{i+1},Y^{\d,\eta_{i+1}}_{\t_{i+1}-}))< h_{j,\eta_{i+1}}(\t_{i+1},Y^{\d,\eta_{i+1}}_{\t_{i+1}-})\le Y^{\d,j}_{\t_{i+1}-}.
 $$
 This, together with (\ref{Y>h11}) and (\ref{Y>h12}), proves (\ref{terminal1}) for
 $j\neq\eta_i,\eta_{i+1}$. $\Box$
 \ms

 Now for each $i$, either (\ref{Ydeta1}) or (\ref{Ydeta2}) is well defined. Therefore, over each $[\t_i,\t_{i+1})$, either (\ref{Ydeta1}) or
  (\ref{Ydeta2}) is wellposed. By applying Corollary \ref{comparison0} and comparison theorem repeatedly, one can easily show
  that:
 \begin{lem}
 \label{upperbound}
 For any admissible strategy $\d$ and any $j$, we have $Y^{\d,j}_t \le Y^j_t$ whenever $Y^{\d,j}_t$ is well defined.
 \qed
 \end{lem}

\subsection{Verification Theorem}
 Moreover, we have :
 \begin{thm}
 \label{verification}
 For $j=1,\cds,m$, we have $\dis Y^{j}_{\l_1} = \esssup_{\d}Y^{\d,{j}}_{\l_1}$.
 \end{thm}

 \noindent {\it Proof.} Fix $\e>0$ and let $D_\e\dfnn \{{i\e}: i=0, 1,\cds\}$. We construct an approximately optimal admissible
  strategy $\d\dfnn \d^\e$ as follows. First, let $\t_0\dfnn \l_1$ and choose $\eta_0$
  such that $A_{\eta_0}\neq\emptyset$. For $i=0,1,\cds$, we define $(\t_{i+1},\eta_{i+1})$ in two cases.
\ms

 \noindent {\it \underline{Case 1}.} If $A_{\eta_i}\neq\emptyset$, set
 $$
 \t_{i+1}\dfnn \inf\{t\ge \t_i: Y^{\eta_i}_t = \max_{k\in A_{\eta_i}} h_{\eta_i,k}(t, Y^k_t)\}\wedge\l_2.
 $$
 If $\t_{i+1}<\l_2$, set $\eta_{i+1}\in A_{\eta_i}$ be the smallest index such that
 \be
 \label{Y=h}
 Y^{\eta_i}_{\t_{i+1}} =  h_{\eta_i,\eta_{i+1}}(\t_{i+1}, Y^{\eta_{i+1}}_{\t_{i+1}}).
 \ee
 Otherwise choose arbitrary $\eta_{i+1}\in A_{\eta_i}$.

\noindent {\it \underline{Case 2}.} If  $A_{\eta_i}=\emptyset$, since
$A_{\eta_0}\neq\emptyset$, we have $i\ge 1$. Set $\eta_{i+1}\dfnn
\eta_{i-1}$. If $\t_i=\l_2$, define $\t_{i+1}\dfnn \l_2$. Now assume
$\t_i<\l_2$. It is more involved to
 define $\t_{i+1}$ in this case. By the definition
  of $\eta_i$, one can check that in this case we must have $A_{\eta_{i-1}}\neq\emptyset$,
  and thus by Case 1, $\eta_i\in A_{\eta_{i-1}}$ and
$$
Y^{\eta_{i-1}}_{\t_i} = h_{\eta_{i-1},\eta_i}(\t_i,
Y^{\eta_i}_{\t_i}).
$$
We claim that, for any $j$ such that $\eta_{i-1}\in A_j$, \be
\label{Y>h} Y^j_{\t_i} > h_{j,\eta_{i-1}}(\t_i,
Y^{\eta_{i-1}}_{\t_i}). \ee In fact, if not, by Assumption
\ref{assump2} (ii), $\eta_i\in A_j\cup\{j\}$ and
$$
Y^j_{\t_i} =  h_{j,\eta_{i-1}}(\t_i, Y^{\eta_{i-1}}_{\t_i}) =
h_{j,\eta_{i-1}}(\t_i, h_{\eta_{i-1},\eta_i}(\t_i,
Y^{\eta_i}_{\t_i})) < h_{j,\eta_i}(\t_i, Y^{\eta_i}_{\t_i}).
$$
This contradicts with (\ref{Yjj}). We now define
$$
\t_{i+1}\dfnn \t^1_{i+1}\wedge \t^2_{i+1} \wedge \l_2;
$$
where $\t^1_{i+1}$ is the smallest number in $D_\e$ such that
$\t^1_{i+1}>\t_i$; and
 $$
 \t^2_{i+1} \dfnn \inf\{t>\t_i: \exists j ~ s.t.~ \eta_{i-1}\in A_j,
 Y^j_{t} =  h_{j,\eta_{i-1}}(t, Y^{\eta_{i-1}}_{t})\}.
 $$

 We claim that, for a.s. $\o$, $\t_n = \l_2$ for $n$ large enough.
 In fact, if $\t_n<\l_2$ for all $n$, let $\dis \t_\infty \dfnn \lim_{n\to\infty}\t_n$.
 In Case 1, (\ref{Y=h}) holds true. In Case 2, if $\t_{i+1} = \t^1_{i+1}$, then $\t_{i+1}\in D_\e$ ;
 and if $\t_{i+1} = \t^2_{i+1}$, then there exists $\hat\eta_{i+1}$ such that $\eta_{i-1}\in A_{\hat\eta_{i+1}}$
 and
 \be
 \label{Y=h2}
 Y^{\hat\eta_{i+1}}_{\t_{i+1}} = h_{\hat\eta_{i+1},\eta_{i-1}}(\t_{i+1},Y^{\eta_{i-1}}_{\t_{i+1}}).
 \ee
 Since $\t_i<\infty$ for all $i$, there can be only finitely many $i$ such that $\t_{i+1}\in D_\e$.
 Therefore, there exists some $n_0$ such that for all $i\ge n_0$, either (\ref{Y=h}) or (\ref{Y=h2})
 holds true. The vector $(\hat \eta_{i+1}, \eta_{i-1}, \eta_i)$ can take only finitely many values,
 then there exist $(j_1,j_2,j_3)$ and an infinite sequence of $i_k$ such that $j_2\in A_{j_1}, j_3\in A_{j_2}$ and
 $$
 \hat \eta_{i_k+1} = j_1,\q \eta_{i_k-1} = j_2,\q \eta_{i_k}=j_3,\q \forall k.
 $$
 By (\ref{Y=h2}) and (\ref{Y=h}) we get
 $$
 Y^{j_1}_{\t_{i_k+1}} = h_{j_1,j_2}(\t_{i_k+1},Y^{j_2}_{\t_{i_k+1}}),\q
 Y^{j_2}_{\t_{i_k}} = h_{j_2,j_3}(\t_{i_k},Y^{j_3}_{\t_{i_k}}),\q \forall k.
 $$
 Send $k\to\infty$, we have
 $$
 Y^{j_1}_{\t_\infty} = h_{j_1,j_2}(\t_\infty,Y^{j_2}_{\t_\infty}),\q
  Y^{j_2}_{\t_\infty} = h_{j_2,j_3}(\t_\infty,Y^{j_3}_{\t_\infty}).
 $$
 Then, by Assumption \ref{assump2} $(ii)$, $j_3\in A_{j_1}\bigcup\{j_1\}$ and
 $$
 Y^{j_1}_{\t_\infty} = h_{j_1,j_2}(\t_\infty, h_{j_2,j_3}(\t_\infty,Y^{j_3}_{\t_\infty}))
 < h_{j_1,j_3}(\t_\infty,Y^{j_3}_{\t_\infty}).
 $$
This contradicts with (\ref{Yjj}). Therefore, $\t_n = \l_2$ for $n$
large enough.

 We now set $\d^{n,\e}\dfnn (\t_0,\cds,\t_n; \eta_0,\cds,\eta_n)$.
  Recall  Definition \ref{admissible}. One can easily check that $\d^{n,\e}$ is an admissible strategy.
Denote
$$
\D Y^j_t \dfnn Y^j_t - Y^{\d^{n,\e},j}_t.
$$
If $i+1=n$, it is obvious that
 \be
 \label{DY0}
 | Y^j_{\t_{i+1}}-Y^{\d^{n,\e},j}_{\t_{i+1}-}|=|\D Y^j_{\t_{i+1}}|.
 \ee

 We now assume $i+1<n$.

 \noindent {\it \underline{Case 1}.} Note that $(Y^j,Z^j,K^j)$ satisfies
\be
 \label{Y1}
 \left\{\ba{lll}
 \dis Y^{j}_t = Y^{j}_{\t_{i+1}}+\int_t^{\t_{i+1}} f_j(s,\overrightarrow{Y}_s,Z^{j}_s)ds- \int_t^{\t_{i+1}} Z^{j}_s dB_s + K^{j}_{\t_{i+1}} - K^{j}_t,~ j\neq \eta_i;\\
 \dis Y^{j}_t \ge \max_{k\in A_j}h_{j,k}(t,Y^{k}_t);\q [Y^{j}_t - \max_{k\in A_j}h_{j,k}(t,Y^{k}_t)]dK^{k}_t = 0,~ j\neq \eta_i;\\
 \dis Y^{\eta_i}_t =  
  Y^{\eta_i}_{\t_{i+1}}+\int_t^{\t_{i+1}}
f_{\eta_i}(s,\overrightarrow{Y}_s,Z^{\eta_i}_s)ds- \int_t^{\t_{i+1}}
Z^{\eta_i}_s dB_s.
 \ea\right.
 \ee
 Compare (\ref{Y1}) and (\ref{Ydeta1}). By induction we have
 \be
 \label{DY11}
 \max_{1\le j\le m} |\D Y^j_{\t_i}|^2 \le  E_{\t_i}\Big\{ e^{C(\t_{i+1}-\t_i)}\max_{1\le j\le m}| Y^j_{\t_{i+1}}-Y^{\d^{n,\e},j}_{\t_{i+1}-}|^2\Big\}.
 \ee
 If $\t_{i+1}=\l_2$, then
 \be
 \label{DY12}
 |Y^j_{\t_{i+1}}-Y^{\d^{n,\e},j}_{\t_{i+1}-}|=|\xi^j_{\l_2}-\xi^j_{\l_2}|=0,\q\forall j.
 \ee
 Assume $\t_{i+1}<\l_2$. Note that $Y^{\d^{n,\e},j}_{\t_{i+1}-}$ is defined by either (\ref{Yi-2}) or (\ref{Yi-1}). In the former case, by (\ref{hLipschitz}) we have
 \beaa
 &&\max_{j\neq \eta_i}| Y^j_{\t_{i+1}}-Y^{\d^{n,\e},j}_{\t_{i+1}-}|=\max_{j\neq\eta_i}|\D Y^j_{\t_{i+1}}|;\\
 && | Y^{\eta_i}_{\t_{i+1}}-Y^{\d^{n,\e},\eta_i}_{\t_{i+1}-}|=| h_{\eta_i,\eta_{i+1}}(\t_{i+1},Y^{\eta_{i+1}}_{\t_{i+1}})- h_{\eta_i,\eta_{i+1}}(\t_{i+1},Y^{\d^{n,\e},\eta_{i+1}}_{\t_{i+1}})|\le |\D Y^{\eta_{i+1}}_{\t_{i+1}}|.
 \eeaa
 Then
 \be
 \label{DY13}
 \max_{1\le j\le m}| Y^j_{\t_{i+1}}-Y^{\d^{n,\e},j}_{\t_{i+1}-}|\le \max_{j\neq\eta_i}| \D Y^j_{\t_{i+1}}|.
 \ee
 In the latter case, recalling Lemma \ref{upperbound} and (\ref{h2}), we have
 \beaa
 &&\max_{j\neq \eta_i, \eta_{i+1}}| Y^j_{\t_{i+1}}-Y^{\d^{n,\e},j}_{\t_{i+1}-}|=\max_{j\neq\eta_i,\eta_{i+1}}|\D Y^j_{\t_{i+1}}|;\\
 &&|Y^{\eta_{i+1}}_{\t_{i+1}}-Y^{\d^{n,\e},\eta_{i+1}}_{\t_{i+1}-}|\le |\D
 Y^{\eta_{i+1}}_{\t_{i+1}}|;\\
 &&|Y^{\eta_{i}}_{\t_{i+1}}-Y^{\d^{n,\e},\eta_{i}}_{\t_{i+1}-}|
 = | h_{\eta_i,\eta_{i+1}}(\t_{i+1},Y^{\eta_{i+1}}_{\t_{i+1}})-
 h_{\eta_i,\eta_{i+1}}(\t_{i+1},Y^{\d^{n,\e},\eta_{i+1}}_{\t_{i+1}-})|\\
&&\qq \le
|Y^{\eta_{i+1}}_{\t_{i+1}}-Y^{\d^{n,\e},\eta_{i+1}}_{\t_{i+1}-}|\le
|\D Y^{\eta_{i+1}}_{\t_{i+1}}|.
 \eeaa
 Thus (\ref{DY13}) also holds true. Therefore, in all the cases we get
 \be
 \label{DY1}
 \max_{1\le j\le m} |\D Y^j_{\t_i}|^2 \le  E_{\t_i}\Big\{ e^{C(\t_{i+1}-\t_i)}\max_{j\neq \eta_i}|\D Y^j_{\t_{i+1}}|^2\Big\}.
 \ee

 \noindent {\it \underline{Case 2}.} Note that  $(Y^j,Z^j,K^j), j\neq \eta_{i-1}$ satisfies
\be
 \label{Y2}
 \left\{\ba{lll}
 \dis Y^{j}_t = Y^{j}_{\t_{i+1}}- \int_t^{\t_{i+1}} Z^{j}_s dB_s + K^{j}_{\t_{i+1}} - K^{j}_t\\
  \dis\qq+\int_t^{\t_{i+1}} \hat f_j(s,Y^{1}_s,\cds,Y^{\eta_{i-1}-1}_s,Y^{\eta_{i-1}+1}_s,\cds, Y^{m}_s,Z^{j}_s)ds;\\
 \dis Y^{j}_t \ge \max_{k\in A_j-\{\eta_{i-1}\}}h_{j,k}(t,Y^{k}_t);\q
 [Y^{j}_t - \max_{k\in A_j-\{\eta_{i-1}\}}h_{j,k}(t,Y^{k}_t)]dK^{k}_t = 0;
 \ea\right.
 \ee
 where
 \bea
 \label{hatf}
 &&\hat f_j(t,y_1,\cds,y_{\eta_{i-1}-1},y_{\eta_{i-1}+1},\cds,
 y_n,z)\\
 &&\q \dfnn \tilde f_j(t,y_1,\cds,y_{\eta_{i-1}-1},y_{\eta_{i-1}+1},\cds,
 y_n,z) + I^j_t;\nonumber\\
 \label{Ijt}
 && I^j_t \dfnn f_j(t,\overrightarrow{Y}_t,Z^j_t)\\
 &&\qq - f_j(t,Y^1_t,\cds,Y^{\eta_{i-1}-1}_t,h_{\eta_{i-1},\eta_i}(\t_i,Y^{\eta_i}_t),Y^{\eta_{i-1}+1}_t,\cds,Y^n_t,Z^j_t)\nonumber.
 \eea
 We note that here $I^j_t$ is considered as a random coefficient. Compare (\ref{Y2}) and (\ref{Ydeta2}).
 Recalling (\ref{Yi-3}), by induction we get
 \be
 \label{DY2}
 \max_{j\neq \eta_{i-1}} |\D Y^j_{\t_i}|^2 \le
 E_{\t_i}\Big\{ e^{C(\t_{i+1}-\t_i)}\max_{j\neq \eta_{i-1}}|\D Y^j_{\t_{i+1}}|^2
 + C\sum_{j\neq\eta_{i-1}} \int_{\t_i}^{\t_{i+1}} |I^j_t| dt\Big\}.
 \ee
 Note that
 $
 Y^{\eta_{i-1}}_{\t_i} = h_{\eta_{i-1},\eta_i}(\t_i,
 Y^{\eta_i}_{\t_i}).
 $
 Then
 \beaa
 |I^j_t|&\le& C\Big|Y^{\eta_{i-1}}_t- h_{\eta_{i-1},\eta_i}(\t_i,Y^{\eta_i}_t)\Big|^2\\
 &\le& C\Big[|Y^{\eta_{i-1}}_t-Y^{\eta_{i-1}}_{\t_i}|^2 + |h_{\eta_{i-1},\eta_i}(\t_i,
 Y^{\eta_i}_{\t_i})-h_{\eta_{i-1},\eta_{i}}(\t_i, Y^{\eta_i}_t)|^2\Big]\\
 &\le& C\Big[|Y^{\eta_{i-1}}_t-Y^{\eta_{i-1}}_{\t_i}|^2 + |Y^{\eta_i}_{\t_i}- Y^{\eta_i}_t|^2\Big]\le
 C\sum_{k=1}^m |Y^k_t - Y^k_{\t_i}|^2.
 \eeaa
 Note that in this case  $\t_{i+1}-\t_i\le \e$. Then
 \be
 \label{Ie}
 |I^j_t| \le C\sum_{k=1}^m\sup_{\l_1\le t_1< t_2\le \l_2: t_2-t_1\le \e}|Y^k_{t_1}-Y^k_{t_2}|^2 \dfnn I_\e.
 \ee
 Thus (\ref{DY2}) implies
 \be
 \label{DY3}
 \max_{j\neq \eta_{i-1}} |\D Y^j_{\t_i}|^2 \le  E_{\t_i}\Big\{ e^{C(\t_{i+1}-\t_i)}\max_{j=1,m}|\D Y^j_{\t_{i+1}}|^2 + I_\e[\t_{i+1}-\t_i]\Big\}.
 \ee

 Now given $A_{\eta_i}\neq\emptyset$, if $A_{\eta_{i+1}}=\emptyset$, by
 (\ref{DY1}) and (\ref{DY3}) we have
 \be
 \label{DY4}
 \max_{1\le j\le m} |\D Y^j_{\t_i}|^2 \le  E_{\t_i}\Big\{ e^{C(\t_{i+2}-\t_i)}\max_{1\le j\le m}|\D Y^j_{\t_{i+2}}|^2 + I_\e[\t_{i+2}-\t_{i+1}]\Big\}.
 \ee
 By Definition \ref{admissible} (v), we have
 $A_{\eta_{i+2}}\neq\emptyset$. Therefore, if $A_{\eta_i}\neq\emptyset$,
 then either $A_{\eta_{i+1}}\neq\emptyset$ and (\ref{DY1}) holds true, or
  $A_{\eta_{i+2}}\neq\emptyset$ and (\ref{DY4}) holds true. Since $A_{\eta_0}\neq \emptyset$, one
  gets immediately that
 $$
 \max_{1\le j\le m} |\D Y^{j}_{\t_0}|^2 \le CE_{\t_0}\Big\{\max_{1\le j\le m}|\D Y^j_{\t_n}|^2+ I_\e\Big\}
 = CE_{\l_1}\Big\{\max_{1\le j\le m}|Y^{0,j}_{\t_n}- Y^j_{\t_n}|^2+ I_\e\Big\} .
 $$
 First send $n\to\infty$. Since $\t_n\to\l_2$, we get
 $$
 Y^{0,j}_{\t_n}\to \xi^j_{\l_2},\q Y^j_{\t_n}\to \xi^j_{\l_2}.
 $$
 By Dominating Convergence Theorem we have
 $$
 \max_{1\le j\le m} |\D Y^{j}_{\l_1}|^2 \le CE_{\l_1}\{I_\e\}.
 $$
 Now send $\e\to 0$. Since $Y^j$ is continuous,  by Dominating Convergence Theorem again we get
 $$
 \lim_{n\to\infty} E_{\l_1}\{I_\e\} = 0.
 $$
 This proves the theorem.
 \qed

\subsection{Proof of Theorem \ref{uniqueness}}

 As mentioned before, we prove the theorem by induction. Assume
 Theorem \ref{uniqueness} holds true for $\mu=m_1-1$. Now assume
 $\mu=m_1$.

 (i) By Theorem \ref{verification}, $Y^j_{\l_1}$ is unique. Similarly $Y^j_t$ is unique for any $t\in [\l_1,\l_2]$.
 By the uniqueness of
the Doob-Meyer decomposition we get $Z^j$ is unique, which further
implies the uniqueness of $K^j$ immediately. \ms

 (ii) For any admissible strategy $\d$, define $\tilde Y^{\d,j}$ similarly
and denote
 $$
 \D Y^{\d,j}_t \dfnn Y^{\d,j}_t - \tilde Y^{\d,j}_t.
 $$
 If $A_{\eta_i}\neq\emptyset$, recalling (\ref{Ydeta1}),  (\ref{Yi-2}), and (\ref{Yi-1}), by induction we
 have:
 $$
 \max_{1\le j\le m} |\D Y^{\d,j}_{\t_i}|^2
 \le E_{\t_i}\Big\{e^{C(\t_{i+1}-\t_i)}\max_{j\neq \eta_i}|\D Y^{\d,j}_{\t_{i+1}}|^2
 +C\int_{\t_i}^{\t_{i+1}}\|\D f_t\|^2 dt\Big\}.
 $$
 If $A_{\eta_i}=\emptyset$, recalling (\ref{Ydeta2}) and (\ref{Yi-3}), by induction we
 have:
 $$
 \max_{j\neq \eta_{i-1}} |\D Y^{\d,j}_{\t_i}|^2\le
  E_{\t_i}\Big\{e^{C(\t_{i+1}-\t_i)}\max_{1\le j\le m}|\D Y^{\d,j}_{\t_{i+1}}|^2
   +C\int_{\t_i}^{\t_{i+1}}\|\D f_t\|^2 dt \Big\}.
 $$
 Put together and note that $A_{\eta_0}\neq\emptyset$, we get:
 $$
 \max_{1\le j\le m} |\D Y^{\d,j}_{\l_1}|^2\le  E_{\l_1}\Big\{e^{C(\l_2-\l_1)}
 \max_{1\le j\le m}|\D \xi^{j}_{\l_2}|^2+C\int_{\l_1}^{\l_2}\|\D f_t\|^2dt\Big\}.
 $$
 Then (ii) follows from Theorem \ref{verification} immediately.
 \qed

\end{document}